\documentclass[12pt]{amsart}
\usepackage{amscd}
\newtheorem{thm}[equation]{Theorem}
\newtheorem{lem}[equation]{Lemma}
\newtheorem{prop}[equation]{Proposition}
\newtheorem{cor}[equation]{Corollary}

\theoremstyle{definition}

\theoremstyle{remark}

\numberwithin{equation}{section}

\newcommand{\pl}{\mathbf{P}^1} % projective line
\newcommand{\pp}{\mathbf{P}^2} % projective plane
\newcommand{\dpp}{\mathbf{\check{P}}^{2}} % dual projective plane %
\newcommand{\hd}{\check{h}} % class of dual line
\newcommand{\ld}{\lambda_d} % class of lift of immersion

\newcommand{\ppmod}[2]{{\overline M}_{#1}(\pp,#2)}
\newcommand{\sempmod}[2]{{\overline M}_{#1}(S_2,#2)}
\newcommand{\lmod}[2]{{\overline M}^2_{#1}(\pp,#2)}
\newcommand{\lmodnd}{\lmod{n}{d}}

\newcommand{\cuspmod}[1]{{\overline C}(\pp,#1)}
\newcommand{\stablefour}{{\overline M}_4}

\newcommand{\ramstack}[2]{{\overline R}_{#1}(#2)}
\newcommand{\tail}{{\overline T}}

\newcommand{\calC}{{\mathcal C}}
\newcommand{\calF}{{\mathcal F}}
\newcommand{\calG}{{\mathcal G}}
\newcommand{\calJ}{{\mathcal J}}
\newcommand{\calK}{{\mathcal K}}
\newcommand{\calM}{{\mathcal M}}
\newcommand{\calN}{{\mathcal N}}
\newcommand{\calO}{{\mathcal O}}
\newcommand{\calP}{{\mathcal P}}
\newcommand{\calR}{{\mathcal R}}
\newcommand{\calT}{{\mathcal T}}
\newcommand{\calX}{{\mathcal X}}

\newcommand{\jrnd}{\calJ{\overline \calR}^2_n(d)}
\newcommand{\jxnt}{\calJ\calX_{n}^2(T)}
\DeclareMathOperator{\cod}{codim}

\begin{document}

\title [Contact Formulas via Stable Maps]
{Contact Formulas for Rational Plane Curves \\ via Stable Maps}

\author[S.\ J.\ Colley]{Susan Jane Colley}
\address{Department of Mathematics,
Oberlin College, Oberlin, Ohio 44074, USA}
\email{sjcolley@math.oberlin.edu}
\thanks{Research of the first author was partly supported 
by National Science Foundation grant \#DMS-9805819.}

\author[L.\ Ernstr{\"o}m]{Lars Ernstr{\"o}m}
\address{Department of Mathematics,
Royal Institute of Technology,
S-100 44 Stockholm, Sweden}
\email{ernstrom@math.kth.se}
\thanks{Research of the second author was partly supported by the Swedish
Natural
Science Research Council.}

\author[G.\ Kennedy]{Gary Kennedy}
\address{Ohio State University at Mansfield, 1680 University Drive,
Mansfield, Ohio 44906, USA}
\email{kennedy@math.ohio-state.edu}
\thanks{The third author was partially supported by a grant from the
Gustafsson Foundation.}


\begin{abstract}
    Extending techniques of \cite{ErnstromKennedy}, we use stable
    maps, and their stable lifts to the Semple bundle variety of
    second-order curvilinear data, to calculate certain
    characteristic numbers for rational plane curves.  These
    characteristic numbers involve first-order (tangency) and
    second-order (inflectional) conditions. Although they may be
    virtual, they may be used as inputs in an enumeratively
    significant formula for the number of rational curves
    having a triple contact with a specified plane curve and passing
    through $3d-3$ general points.
\end{abstract}

\maketitle

%%%% **** The text of the paper starts here **** %%%% %
%
%
%
\section{Introduction}
\label{intro}

\par
Two plane curves are said to have a \emph{triple contact} at a point if
both curves are smooth at the point and if their intersection number
there is 3. (For example, a curve has a triple contact with its
tangent line at each flex.) For a specified curve $C$, one may ask
how many rational plane curves of degree~$d$ have a triple contact
with $C$ and also pass through $3d-3$ general points of the plane. Our
formula~(\ref{final}) shows that this \emph{contact number} is the inner
product of
two vectors, the first of which consists of the degree~$c$ of the
specified curve, its class $\check{c}$, and the number of
cusps $\kappa$. The second vector consists of three
\emph{second-order characteristic numbers} of the family of rational plane
curves of degree~$d$.
Thus to give a satisfactory answer to the question---and to questions
of a similar nature---one needs a method
of determining these characteristic numbers.
In this paper we show how to compute (in every degree)
thirteen characteristic numbers,
including the three needed in the contact formula, by developing a
recursive procedure for calculating the \emph{Gromov-Witten
invariants} of certain stacks of stable maps, and by showing that
these invariants agree with the characteristic numbers.

\par
In {\S}\ref{semple} we
briefly discuss the Semple bundle construction of a complete variety
of second-order curvilinear data of $\pp$.  In {\S}\ref{stable} we
introduce the notion of second-order stable lifts of maps to $\pp$ and
their associated stacks, and present three key examples.
We introduce the relevant Gromov-Witten invariants in
{\S}\ref{GWinv}, and prove that
they agree with the characteristic numbers of \cite{CKsimult}.
In {\S}\ref{sbd} we study the special boundary divisors on the stack
of second-order stable lifts and describe in detail those components
that actually matter in the subsequent computations.
In {\S}\ref{potpde} we construct various generating functions for
Gromov-Witten invariants, called potentials,
and show how
they are related by a basic differential equation. Using this equation
in {\S}\ref{recform}, we describe a recursive process which
determines, in every degree, the thirteen Gromov-Witten invariants;
we give a complete table of values through degree~$6$. Finally, in
{\S}\ref{contactform}, we present two ``contact formulas''
(generalizing Theorem 4 of \cite{CKsimult}) in which the Gromov-Witten
invariants are the inputs, and prove that the outputs are
enumeratively significant. We also demonstrate the second formula by using it
to calculate the answer to the question posed at the beginning.

\par To obtain explicit answers to related questions,
we would need to know a larger set of
characteristic numbers.
Our recursion determines, in each degree, those thirteen characteristic
numbers specified by at
least $3d-3$ ``point conditions,''
beginning (as input) with those in degree $1$. For plane conics,
there are $153$ second-order
characteristic numbers which are nontrivial, i.e., specified
by a selection of conditions including neither divisors nor the
fundamental class.
If we knew the full set of values (and one could certainly hope
to find them by \textit{ad hoc} methods), our recursive scheme
would then compute the corresponding $153$ numbers in every
degree, and our
Theorem \ref{mixedcontactprop} would then yield an explicit formula for
the number of rational plane curves having a triple contact with each
of two
specified plane curves and passing through $3d-5$ points.
Another possible generalization would be to discover and validate
formulas for higher-order contact. Semple's construction actually
creates a tower of $\pl$-bundles
parametrizing higher-order curvilinear data, which could presumably
be used in such a project. (See \cite{CKtrip} or \cite{CKsimult}.)
There seem to be, however, two
obstructions. First, we do not know how to prove the analogue of
Proposition \ref{irrellemma}, the essential tool in classifying the
relevant components of the special boundary divisors on the stack of
stable lifts.
Second, the number of $PGL(2)$-orbits on the Semple bundle variety
grows larger with the order of data; this makes it difficult to use
transversality theory to establish enumerative significance.
(In fact since the dimension of the variety of $n$th-order data is
$n+2$ there are infinitely many orbits when  $n>6$.) For 
similar reasons, it will probably be difficult to extend our methods to
the study of higher-dimensional varieties or even to nonrational surfaces;
even in the case of other rational surfaces we anticipate complications.
\par
Several people have suggested the possibility of using gravitational
descendants to deal with problems of higher-order contact.
In discussing this approach, however, Gathmann shows by an example
\cite[p.\ 27]{gathmann}
that it will have to deal with spurious contributions
arising from singular curves. We believe that the stable lift method
avoids this difficulty.

\par
We thank the following people for their interest and assistance:
David Cox, Andreas Gathmann, Morten Honsen, Alexandre Kabanov, Lee McEwan,
Rahul Pandharipande, Israel Vain\-sencher, and Ravi Vakil.  We are
also grateful to the referee for providing many valuable suggestions.

\section{The variety of second-order curvilinear data}
\label{semple}

\par
We provide a quick summary of the salient features of $S_2$, the variety of
second-order curvilinear data of $\pp$, introduced
by Semple in \cite{Semple}. For additional details, see \cite{CKtrip} or
\cite{CKsimult}.

\par

Let $S_1$ be the total space of the projectivized tangent bundle of
$\pp$; equivalently, $S_1$ is the incidence correspondence of points and
lines. Let
$f_{1}\colon S_1 \to \pp$ be the tautological projection, so that a fiber
represents the pencil of tangent directions at a point of the plane.
The \emph{focal plane} at
a point $x \in S_1$ is defined to be the preimage, via the derivative
$df_1$, of the line in $T_{f_{1}(x)}\pp$ represented by $x$.  This
local construction gives rise to a global rank 2 \emph{focal plane
bundle} $\calF_1$ and we define $S_2$, the \emph{Semple bundle variety
of second-order data of $\pp$}, to be $\mathbf{P}\calF_{1}$, the
total space of the projectivized bundle.  The
Semple bundle variety $S_2$ is equipped with a tautological projection
$f_2\colon S_2 \to S_1$. Note that $S_2$ is a subvariety of $\mathbf{P}T(S_1)$,
the projectivized tangent bundle of $S_1$.
The relative tangent bundle $T(S_1/\pp)$ is naturally a rank $1$ subbundle
of $\calF_1$, so that $\mathbf{P}T(S_1/\pp)$ is naturally a section of $S_2$,
which we call the {\em divisor at infinity} and denote by $I$.

\par
If $C$ is a reduced curve in $\pp$, there is a rational map from $C$ to
$S_1$ that sends a
nonsingular point $x$ of $C$ to the point of $S_1$ representing the tangent
direction to $C$ at $x$.  We call the closure of the image of this map the
\emph{lift}
of $C$ and denote it by $C_1$. Similarly, $C_1$ may be lifted to a
curve $C_2$ in $\mathbf{P}T(S_1)$. It is not
difficult to see that the tangent direction to $C_1$ at a
nonsingular point must be in the focal plane.
Consequently, the second-order lift $C_2$ is a curve not just in
$\mathbf{P}T(S_1)$, but in $S_2$.
If $C_2$ passes through a point $x$ of $S_2$, we sometimes say
that the germ of $C$ \emph{represents} $x$. For example, the germ at
the origin of $y^2=x^3$ represents a point on the divisor at
infinity. Intuitively, points on this divisor are represented by
curve germs with ``infinite curvature.''

\par Families of curves may also be lifted. If $\calC
\subset \pp \times T$ is a family of plane curves with reduced
general member, then we lift the reduced members of $\calC$ to $S_2$,
obtaining a rational lifting map from $\calC$ to $S_2 \times T$.  We
define $\calC_2 \subset S_2 \times T$ to be the closure over $\calC$
of the image of this lifting map.

\par
In the construction of the Semple bundle variety we call $\pp$ the
\emph{base}. One may also construct a Semple bundle variety using the dual
projective plane $\dpp$ as base. Note that $\mathbf{P}T(\pp)$
and $\mathbf{P}T(\dpp)$ are both isomorphic to the point-line incidence
correspondence. Furthermore, if $(x,l)$ is a point of the incidence
correspondence then the preimages of $x \in \pp$ and $l \in \dpp$
are curves through $(x,l)$ and the focal plane at $(x,l)$ is the plane
spanned by the tangent lines to these curves. Hence in both
constructions the focal plane bundle is the same, and thus the two
second-order Semple bundle varieties are the same. There are, however,
two disjoint divisors at infinity. We shall continue to denote by $I$ the
divisor at infinity that arises from using $\pp$ as base, but denote
by $Z$ the divisor at infinity that arises from using $\dpp$ as base.
Note that the lift of a smooth curve meets $Z$ above $x \in \pp$ if
and only if $x$ is a flex.
\par
Since $S_2 \to S_1 \to \pp$ is a sequence of $\pl$-bundles, we may
regard  $A^*(\pp)$ and $A^*(S_1)$ as subrings of the intersection
ring $A^*(S_2)$.
Let $h$ and $\hd$ denote the hyperplane classes on
$\pp$ and $\dpp$. Let $i$ denote the
class of the divisor at infinity $I$, and $z$ the class of $Z$.
Then, according to {\S}3 of \cite{CKtrip}, the Chow ring $A^*(S_2)$ is
generated by $h$, $\hd$, and $i$, subject to the relations
$$
h^{3} = \hd^{3} = h^{2} - h\hd + \hd^{2} = 0, \quad i^2 = 3(h-\hd)i,
$$
and the class $z$ satisfies these equations:
\begin{equation*}
i - z = 3(h-\hd), \quad iz = 0.
\end{equation*}

\par

\section{Stable lifts and their stacks}
\label{stable}

\par
We assume that the reader is familiar with the basic notions of stable
maps \cite{FultonP}. Throughout this paper we will consider only stable
maps from
curves of arithmetic genus zero; thus we will omit the subscript $0$ from the
standard notation for moduli stacks of stable maps. Given a stable map
from a tree of $\pl$'s, we will use the term
\emph{twig} to refer to one of the $\pl$'s of the tree, reserving the
term \emph{component} to refer to an irreducible component of an
appropriate moduli space or to a divisor on a moduli space in order
to avoid confusion. We will call the nodes of the source curve
\emph{attachment points}.
\par
If $\mu\colon (\pl, p_1,\ldots,p_n) \to \pp$ is a nonconstant map, then we
define its \emph{strict lift} $\mu_1^{\text{strict}}\colon (\pl,
p_1,\ldots,p_n) \to S_1$ by associating to each point $x$ the
direction of $d\mu(x)$ at $\mu(x)$. Note that this map
is defined even at singular points of $\mu$.
Repeating the
construction, we obtain a map  $\mu_2^{\text{strict}}\colon (\pl,
p_1,\ldots,p_n) \to S_2$, called the
\emph{second-order strict lift} of $\mu$.
(One needs to verify that the image of this map actually lies on $S_2$
inside of the projectivized tangent bundle $\mathbf{P}T(S_1)$.)
If $\mu$ is an immersion, then its degree~$d$ determines the homology class
$\ld$ of $\mu_2^{\text{strict}}(\pl)$.
Thus we have a rational \emph{lifting
map} from the moduli stack $\ppmod{n}{d}$ to $\sempmod{n}{\ld}$, and
we define the \emph{stack of second-order stable lifts}, denoted
$\lmodnd$, to be the closure of its image.  The lifting map is
birational onto $\lmodnd$, and its inverse is inclusion followed
by projection.  Given a point $\mu$ of $\ppmod{n}{d}$, we will call
the points of its preimage in $\lmodnd$ its \emph{second-order
stable lifts}.

\par For an immersion, the only possible stable lift is 
the strict lift. For other sorts of stable maps, however, 
the strict lift will be only
a subset of the twigs of a stable lift; for examples, see Propositions
\ref{node}, \ref{cusp}, and \ref{doublecover} below. Furthermore, some
stable maps have more than
one stable lift. This happens, for example, with the triple cover of a
line in $\pp$ ramified at just two points; we will, however, make no
use of this phenomenon.
\par
Suppose that $\nu\colon \pl\to S_1$ is an isomorphism onto a fiber of $S_1$
over a point of $\pp$. Then the strict lift from $\pl$ to
$\mathbf{P}T(S_1)$ actually lands inside $S_2$. Indeed, since the
composite $f_1 \circ \nu$ is a constant map, the image of the strict
lift is even inside the divisor at infinity $I$. As we will see
below, such strict lifts occur as twigs of certain stable lifts.
\par
\begin{prop}\label{fibersquare}
Suppose that $n \geq m$. Then the following diagram (in which the
horizontal morphism are inclusions and the vertical morphisms are
forgetful) is a fiber square:
$$
\begin{CD}
\lmodnd @>>> \sempmod{n}{\ld} \\
@VVV        @VVV\\
\lmod{m}{d} @>>> \sempmod{m}{\ld}
\end{CD}
$$
\end{prop}
\begin{proof}
Over a point of $\lmod{m}{d}$ representing the strict lift of
an $m$-pointed immersion $\pl \to \pp$, the substack of $\lmodnd$
representing the same immersion together with an additional $m-n$
markings, is dense in the fiber of $\sempmod{n}{\ld}$.
\end{proof}

\par The divisors on $\ppmod{n}{d}$ representing maps
from reducible curves are called \emph{boundary divisors}.  There is
one such divisor $D(A_1,A_2;d_1,d_2)$ for each partition
$\{1,\dots,n\}=A_1
\cup A_2$ and each partition $d=d_1+d_2$.  (If $d_1=0$ then $A_1$ must
have at least two elements; if $d_2=0$ then $A_2$ must have at least two
elements.) A general point of $D(A_1,A_2;d_1,d_2)$ represents a map $\mu\colon
C \to \pp$ from a two-twig curve; the first twig carries the markings
indexed by $A_1$, and its image in $\pp$ is a curve of degree~$d_1$; the
second twig carries the markings indexed by $A_2$, and its image in
$\pp$ is a curve of degree~$d_2$.
(See \cite[p.\ 51]{FultonP} for a drawing.)
Let $p \in \pp$ be the image of the
attachment point.
\par

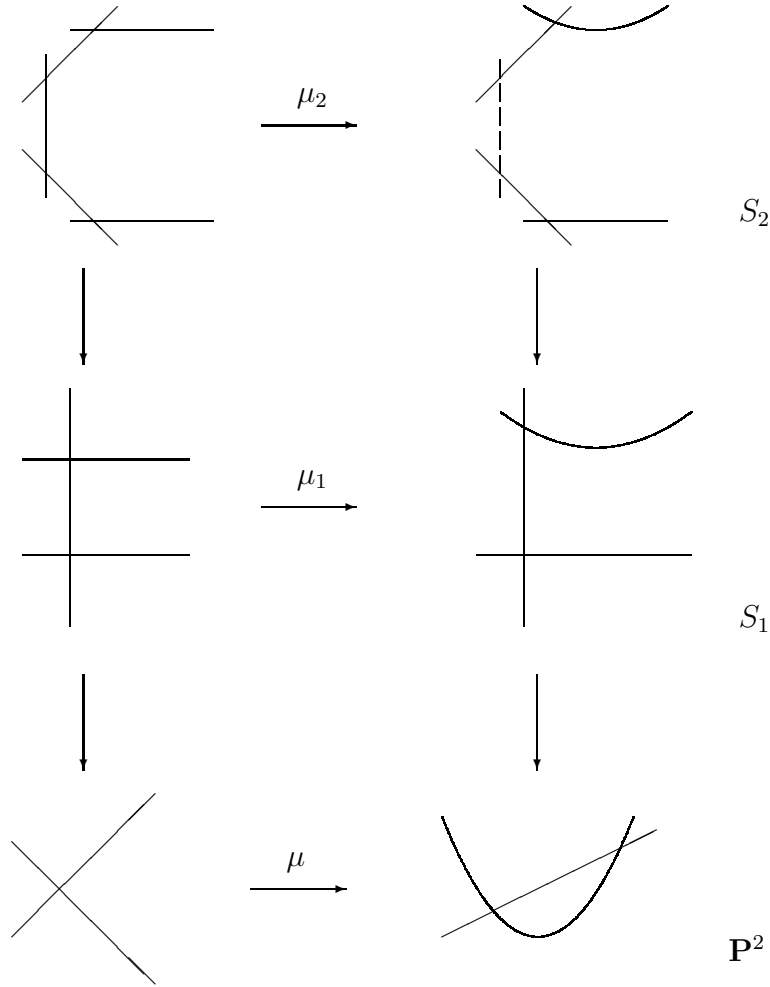
\begin{figure}
   \setlength{\unitlength}{0.125in}
   \begin{center}
   \begin{picture}(32,42)
      \put(3,31){\vector(0,-1){4}}
      \put(3,14){\vector(0,-1){4}}
      \put(22,31){\vector(0,-1){4}}
      \put(22,14){\vector(0,-1){4}}
      \put(0,3){\line(1,1){6}}
      \put(0,7){\line(1,-1){6}}
      \put(10,5){\vector(1,0){4}}
      \put(11.5,6){$\mu$}
      \put(18,3){\line(2,1){9}}
      \qbezier(18,8)(22,-2)(26,8)
      \put(30,2){$\pp$}
   \put(0,16){
      \begin{picture}(32,10)
      \put(2,0){\line(0,1){10}}
      \put(0,3){\line(1,0){7}}
      \put(0,7){\line(1,0){7}}
      \put(10,5){\vector(1,0){4}}
      \put(11.5,6){$\mu_1$}
      \put(21,0){\line(0,1){10}}
      \put(19,3){\line(1,0){9}}
      \qbezier(20,9)(24,6)(28,9)
      \put(30,0){$S_1$}
      \end{picture}}
   \put(0,33){
      \begin{picture}(32,9)
      \put(2,0){\line(1,0){6}}
      \put(2,8){\line(1,0){6}}
      \put(0,3){\line(1,-1){4}}
      \put(0,5){\line(1,1){4}}
      \put(1,1){\line(0,1){6}}
      \put(10,4){\vector(1,0){4}}
      \put(11.5,5){$\mu_2$}
      \put(21,0){\line(1,0){6}}
      \put(19,3){\line(1,-1){4}}
      \put(19,5){\line(1,1){4}}
      \multiput(20,1)(0,1){6}{\line(0,1){.75}}
      \qbezier(21,9)(24,7)(27,9)
      \put(30,0){$S_2$}
      \end{picture}}
   \end{picture}
   \end{center}
\caption{The second-order stable lift $\mu_2$ of the map $\mu$ in
Proposition
\ref{node}.  The vertical twig of the source of $\mu_2$ maps into the
divisor at infinity in $S_2$; its image is shown dashed.}
\label{nodepic}
\end{figure}

\begin{prop}\label{node}
Suppose that $d_1$ and $d_2$ are both positive. Then a general
element $\mu$ of $D(\emptyset,\emptyset;d_1,d_2)$ has a unique
second-order stable lift $\mu_2$. The source
curve is a chain of five twigs. On each peripheral twig, $\mu_2$ is the
second-order strict lift of one of the two twigs of $\mu$.
On the central twig, $\mu_2$ is a map of degree two into the divisor at
infinity $I$ and onto the pencil of tangent directions at $p$, with the
points of attachment mapping to the tangent directions of the two
image curves; the map is ramified at these attachment points. On
each of the remaining two twigs, $\mu_2$ is a map of degree three
onto the fiber of $S_2$ over the point of $S_1$ representing
the tangent direction at $p$ of one of the two image curves. It is
ramified only at the points of attachment. (See Figure \ref{nodepic}.)
\end{prop}
%This is the node figure
\begin{proof}
Since $\ppmod{0}{d}$ is smooth and birationally equivalent to
$\lmod{0}{d}$, the general point of the divisor
$D(\emptyset,\emptyset;d_1,d_2)$ has a
unique preimage.  To see the nature of the stable lift, we follow the
semistable reduction recipe of Fulton and
Pandharipande \cite[pp.\ 64--65]{FultonP}.
We may assume that the source
curve of $\mu$ has just two twigs, that on each twig the map is an
immersion, and that the maps are transverse at the attachment point.
We may create a family in which the source curves are the plane
curves $xy=\epsilon$ and in which the central member is $\mu$. 
Away from the central member, the strict lifts piece together to
give a map to $S_1$, and we now proceed to compute the
``limiting member'' over $\epsilon=0$. (Since the moduli space
is separated and proper, there is a unique such limit.)
In fact the map to $S_1$ extends to
every point of $xy=0$ except the point of attachment, and we
can remove the indeterminacy by one blowup.  
The central member will
now have three twigs, and a calculation shows that on the central twig
the map to $S_1$ is an isomorphism onto a fiber over a point of $\pp$,
while on each peripheral twig it is the strict lift of the original
map. 
The calculation also shows that the map to $S_2$ is
indeterminate at the two attachment points.  Again the indeterminacy
is removed by blowing up the two points, and now the central member
has five twigs arranged in a chain.  The map to $S_2$ takes the newly
introduced twigs onto fibers over points over $S_1$.  On the central
twig the map is the lift of the map at the previous stage; thus it is
a map into the divisor at infinity.  On each peripheral twig the map
is the second-order strict lift of the original map.
\par
The central
member is not reduced, however: its central twig has multiplicity 2
and the adjacent twigs have multiplicity 3.  Thus, following the
recipe, we introduce the base change $\epsilon=\delta^6$ and then
normalize the resulting surface.  After these operations we find that
the central member is still a chain of five twigs, that the map on the
central member is a degree 2 cover, that the maps on the adjacent
twigs are degree 3 covers, and that all ramification is concentrated
at the points of attachment.
\end{proof}

Now consider maps of degree~$d \geq 3$ from $\pl$ to $\pp$
which are immersions except at one point,
at which the map onto the image curve is ramified.
Let $\cuspmod{d}$ be the closure of
the locus on $\ppmod{0}{d}$ representing such maps; it is a
divisor.  For such a map, let $p \in \pp$ indicate the image
of the ramification point, and let $q \in S_1$ indicate its
image under the first-order strict lift (i.e, the point of
$S_1$ representing the cusp tangent direction at $p$).

\begin{prop}\label{cusp}
A general element $\mu$ of $\cuspmod{d}$ has a unique second-order stable
lift $\mu_2$. The source curve has five twigs. On one of these twigs,
which is attached to each of the others, the map to $S_2$ is the
constant map to the point of $I$ representing the tangent directions
to $\pp$ at $p$. The other four twig maps are

\begin{itemize}

\item the strict lift of $\mu$ (attached at a point labeled $a$)

\item a map of degree one into the divisor at
infinity $I$ and onto the pencil of tangent directions at $p$ (attached at
a point $b$)

\item two maps of degree one onto the fiber of $S_2$ over $q$  (attached at
points $c$ and $d$)

\end{itemize}
The attachment points are arranged in such a way that there could
exist a double cover of $\pl$ ramified at $a$ and $b$, with $c$ and $d$ as
the points
of a fiber. (In other words, the cross-ratio $[a,b,c,d]$ is 2.
See Figure \ref{cusppic}.)
\end{prop}
%This is the cusp figure
\begin{figure}
   \begin{center}
   \setlength{\unitlength}{0.125in}
      \begin{picture}(28,42)
      \put(0,5){\line(1,0){8}}
      \put(10,5){\vector(1,0){4}}
      \put(11.5,6){$\mu$}
      \qbezier(21,0)(21,5)(18,5)
      \qbezier(18,5)(21,5)(21,10)
      \put(4,27){\vector(0,-1){4}}
      \put(4,12){\vector(0,-1){4}}
      \put(18,27){\vector(0,-1){4}}
      \put(18,12){\vector(0,-1){4}}
      \put(26,0){$\pp$}
   \put(0,13){
      \begin{picture}(28,8)
      \put(4,0){\line(0,1){8}}
      \put(0,4){\line(1,0){8}}
      \put(10,4){\vector(1,0){4}}
      \put(11.5,5){$\mu_1$}
      \put(18,0){\line(0,1){8}}
      \qbezier(23,7)(18,6)(18,4)
      \qbezier(18,4)(18,2)(23,1)
      \put(26,0){$S_1$}
      \end{picture}}
  \put(0,28){
      \begin{picture}(28,12)
      \put(0,2){\line(1,0){8}}
      \put(0,4){\line(1,0){8}}
      \put(0,6){\line(1,0){8}}
      \put(0,8){\line(1,0){8}}
      \put(2.75,8.25){$a$}
      \put(2.75,6.25){$b$}
      \put(2.75,4.25){$c$}
      \put(2.75,2.25){$d$}
      \put(3,10.5){contracts}
      \put(10,4){\vector(1,0){4}}
      \put(11.5,5){$\mu_2$}
      \multiput(18,0)(0,1){10}{\line(0,1){.75}}
      \put(15,1){\line(1,1){7}}
      \qbezier(14,7)(15,5)(18,4)
      \qbezier(18,4)(21,3)(23,5)
      \thicklines
      \put(4,0){\line(0,1){10}}
      \put(26,0){$S_2$}
      \end{picture}}
      \end{picture}
  \end{center}
\caption{The second-order stable lift $\mu_2$ of the map $\mu$ in
Proposition \ref{cusp}.  The vertical twig in the source of $\mu_2$
maps to a point of $S_2$.  The horizontal twig through $b$ maps into
the divisor at infinity; its image is shown dashed.  The horizontal
twigs through $c$ and $d$ both map onto the fiber of $S_2$ over the
cusp tangent direction.}
\label{cusppic}
\end{figure}
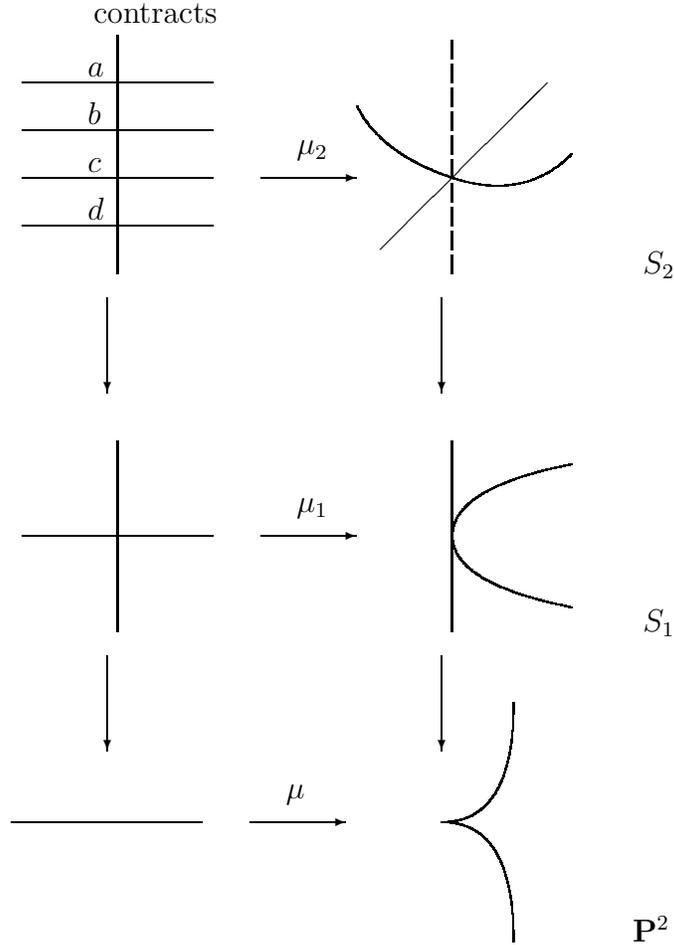
\begin{proof}
Again we know that the general point of the divisor $\cuspmod{d}$ has a
unique preimage, and again we see the nature of the stable lift by
following the recipe of Fulton and Pandharipande. This time we may
begin with a family whose source is the product of $\pl$ and a curve,
with the map to $\pp$ being an immersion except on the special
member. The maps to $S_1$ and $S_2$ obtained by piecing together the
strict lifts are defined everywhere but at the ramification point of
the special member. At this point and at $p$ we may choose local
coordinates so that the special member is $t \mapsto (t^2,t^3)$, then
extend to the family of maps
$$
(\epsilon,t) \mapsto (t^2,t^3-\epsilon t).
$$
The indeterminacy in the map to $S_1$ is resolved by one
blowup. On one of the two twigs of the new special member we
have the strict lift; on the other we have an isomorphism
onto the fiber of $S_1$ over $p$; the attachment point maps
to $q$; and the image curves are tangent there. To resolve
the indeterminacy in the map to $S_2$, we first blow up the
new attachment point. Then the new special member has a
chain of three twigs, with the map on the central one being
constant. There is still, however, indeterminacy at a point
on this twig. One more blowup creates a fourth twig, mapping
onto the fiber of $S_2$ over $q$. Both this twig and the
central one occur in the special member with multiplicity
$2$. Thus when we make the base change $\epsilon=\delta^2$
and then normalize the resulting surface we obtain two
inverse images of the fourth twig, and on the central twig
the normalization map is a double cover ramified at the
points of attachment to the other two twigs.
\par

\end{proof}
\par
The case $d=2$ is special. If $\mu\colon \pl\to\pp$ is not an immersion,
then it must be a double cover of a line. For such a map, let $a$ and
$b$ denote the branch points in $\pp$. Let $\cuspmod{2}$ be the closure
of the locus on $\ppmod{0}{d}$ representing such maps; it is a divisor.

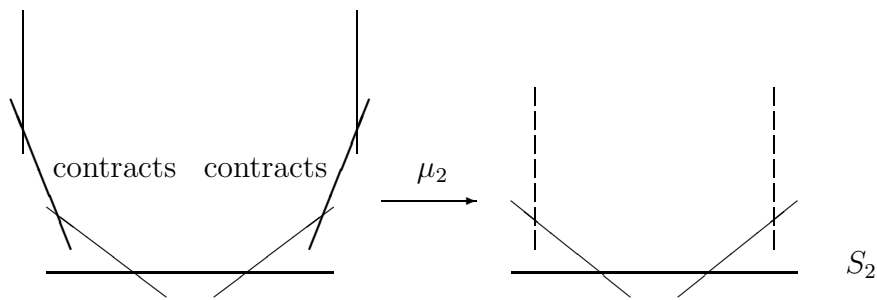
\begin{figure}
   \begin{center}
   \setlength{\unitlength}{0.125in}
      \begin{picture}(38,12)
        %left part
        \put(8,0){\line(-4,3){5}}
        \put(2,6){\line(0,1){6}}
        \put(10,0){\line(4,3){5}}
        \put(16,6){\line(0,1){6}}
        \put(3,1){\line(1,0){12}}
        \put(3.25,5){\makebox(5.5,3.5)[bl]{contracts}}
	\put(9.25,5){\makebox(5.5,3.5)[br]{contracts}}
        \put(17,4){\vector(1,0){4}}
        \put(18.5,5){$\mu_2$}
        \thicklines
        \put(4,2){\line(-2,5){2.5}}
        \put(14,2){\line(2,5){2.5}}
        %right part
        \put(22,0){
          \begin{picture}(16,9)
          \thinlines
          \put(5,0){\line(-5,4){5}}
          \put(7,0){\line(5,4){5}}
          \multiput(1,2)(0,1){7}{\line(0,1){.75}}
          \multiput(11,2)(0,1){7}{\line(0,1){.75}}
          \put(0,1){\line(1,0){12}}
          \put(14,1){$S_2$}
          \end{picture}}
       \end{picture}
   \end{center}
   \caption{The second-order stable lift of a double cover of a line
   in $\pp$.  Each of the two twigs of attached to the peripheral
   twigs is mapped to a point of $S_2$.  The two peripheral twigs map
   into the divisor at infinity; their images are shown dashed.}
   \label{doublepic}
   \end{figure}

\begin{prop}\label{doublecover}
A double cover of a line has a unique stable lift. The source curve
is a chain of seven twigs. On the central twig the map is the strict
lift of the double cover. On each of the adjacent twigs the map is a triple
cover of a fiber of $S_2$ over $S_1$, ramified only at the attachment
points. On
each peripheral twig we have a map into the divisor at infinity and
onto the pencil of directions at $a$ or $b$. On each of the remaining
twigs the map is constant. (See Figure \ref{doublepic}.)
\end{prop}
%This is the figure for the stable lift of the double cover
\begin{proof}
Again we know that the general point of the divisor $\cuspmod{2}$ has a
unique preimage, and again we see the nature of the stable lift by
following the recipe. We may choose local coordinates so that the
family of maps to $\pl$ is
$$
(\epsilon,t) \mapsto (t^2,\epsilon t).
$$
Over the origin, the indeterminacy of the map to $S_1$ is
resolved by one blowup. The map to $S_2$ is still indeterminate at
the new attachment point, and two further blowups are required. (The
first blowup creates a twig on which the map is constant, and the second
one creates a twig mapping to a fiber of $S_2$.) A
similar sequence of blowups over the other branch point creates a
chain or seven twigs. Two of these twigs (those on which the map is
constant) occur in the special member
with multiplicity 2, and two others (those created by the last
blowups) occur with multiplicity 3. Thus we make
the base change $\epsilon=\delta^6$ and then normalize, creating
double and triple covers with ramification concentrated at the points
of attachment.
\end{proof}
\par

\section{Gromov-Witten invariants and characteristic numbers}
\label{GWinv}
Here we define the Gromov-Witten invariants and characteristic
numbers of the family of rational plane curves of degree~$d$, and
show that they agree.
\par
As before, let $\lmodnd$ be the stack of second-order stable
lifts.
Let $e_k$ denote evaluation at the $k$th point.
Suppose that $\gamma_1, \dots, \gamma_n$ are classes in
$A^*(S_2)$. Then the associated {\em Gromov-Witten invariant} is the
number
$$
\left\langle \prod_{k=1}^n \gamma_k \right\rangle_d=\int \prod_{k=1}^n
e_k^*(\gamma_k)\cap [\lmodnd].
$$
Here are two basic properties:
\begin{enumerate}
 \item If $\gamma_n=1$ then $\left\langle \prod\gamma_k\right\rangle_d=0$.
 \item If $\gamma_n$ is a divisor class, then
$$
\left\langle  \prod_{k=1}^n\gamma_k\right\rangle_d=(\int \gamma_n\cap
\ld)\left\langle
\prod_{k=1}^{n-1}\gamma_k\right\rangle_d$$ where $\ld \in A_1(S_2)$ is
the class of the image of the second-order
strict lift of an immersion of degree~$d$.
In particular,
\begin{equation*}
\left\langle  \prod_{k=1}^n\gamma_k\right\rangle_d =
  \begin{cases}
   d\left\langle  \prod_{k=1}^{n-1}\gamma_k\right\rangle_d
        &\text{if $\gamma_n = h$},  \\
   (2d-2)\left\langle  \prod_{k=1}^{n-1}\gamma_k\right\rangle_d
        &\text{if $\gamma_n = \hd$},  \\
   0    &\text{if $\gamma_n = i$}, \\
   (3d-6)\left\langle  \prod_{k=1}^{n-1}\gamma_k\right\rangle_d
        &\text{if $\gamma_n = z$}.
    \end{cases}
\end{equation*}
(The lift of an immersion misses the divisor at infinity; the last
formula is a consequence of the others and the identity $z=i-3h-3\hd$.)
\end{enumerate}
\par
Now let $\calR(d) \subset \mathbf{P}^N \times \pp$
denote the total space of the
family of rational plane curves of degree~$d$, over the
parameter space $T(d)\subset \mathbf{P}^N$, where
$N= \binom{d+2}{2}-1$.
Let ${\overline \calR}(d)$ and ${\overline T}(d)$ denote the closures.
Let
$$
{\overline \calR}^2_n(d) =
({\overline \calR}(d))_2
\times_{{\overline T}(d)}({\overline \calR}(d))_2
\times_{{\overline T}(d)}
\cdots \times_{{\overline T}(d)}
({\overline \calR}(d))_2,
$$
the $n$-fold fiber product of the second-order lift of this
family.
It is a subvariety of $S_2 \times S_2 \times \cdots \times
S_2 \times
{\overline T}(d)$. A general point of ${\overline T}(d)$
represents a curve
with a unique second-order lift, and the fiber of
${\overline \calR}^2_n(d)$
over this point is a variety of dimension $n$ representing
$n$-tuples of
points on the lift. Thus there is a unique component of
${\overline
\calR}^2_n(d)$ dominating
$$
{\overline \calR}_n(d) = {\overline \calR}(d)
\times_{{\overline T}(d)}{\overline \calR}(d)
\times_{{\overline T}(d)}
\cdots
\times_{{\overline T}(d)} {\overline \calR}(d),
$$
which we call the \emph{join}, and denote by $\jrnd$.
(There may be other components. For example, a double line
has a
$2$-parameter family of second-order lifts.
Thus if $d=2$ and $n \geq 3$ the fiber product has another
component
supported over the locus of double lines; if $n \geq 4$ this
component even has a larger dimension than the join.)
\par
Let $e_k\colon \jrnd \to S_2$ be projection onto the $k$th
factor.  Then we define the \emph{characteristic numbers} by
$$
\left\{\prod_{k=1}^n \gamma_k \right\}_d = \int
\prod_{k=1}^n
e_k^*(\gamma_k)\cap [\jrnd].
$$
These characteristic numbers need not have enumerative
significance. For example,
$$
\left\{ h^2 \cdot h^2 \cdot z \right\}_1 = -3,
$$
as one can see by an excess intersection calculation or by
replacing
$z$ by $i-3h+3\hd$.
Nevertheless these numbers are the appropriate inputs for
the contact
formulas of {\S}\ref{contactform}, whose output does have
enumerative significance.
\par

\begin{thm}\label{comparison} For each $d$, $n$,
    and classes $\gamma_1, \dots, \gamma_n$,
$$
\left\langle \prod_{k=1}^n \gamma_k \right\rangle_d=\left\{\prod_{k=1}^n
\gamma_k \right\}_d.
$$
\end{thm}

\begin{proof}
There is a birational morphism from $\lmodnd$ to $\jrnd$ which associates
to the
stable lift of an immersion its image curve in $S_2$
together with the images of the $n$ markings. Denote the graph of this
morphism by $\calG$. Then two applications of the projection formula
show that the Gromov-Witten invariant and the characteristic number
both equal
$$
\int \prod_{k=1}^n e_k^*(\gamma_k)\cap [\calG].
$$
\end{proof}

\section{Special boundary divisors}
\label{sbd}

Suppose that $n \geq 4$. Consider the ``forgetful'' morphism
from $\sempmod{n}{\ld}$ to
$\stablefour$, the moduli space of 4-pointed stable genus $0$ curves, which
associates to a stable map its source curve with
only the first four markings retained. (Unstable twigs are contracted
as necessary.) Composing with inclusion, we have a morphism $\lmod{n}{d} \to
\stablefour$. The space $\stablefour$ is isomorphic to $\pl$. It has a
distinguished point $P(12\mid 34)$ representing the two-twig curve
having the first two markings on one twig and the latter two on the
other; similarly there are two other distinguished points $P(13\mid 24)$ and
$P(14\mid 23)$. Thus on $\lmod{n}{d}$ there are three linearly
equivalent divisors $D(12\mid 34)$, $D(13\mid 24)$, and $D(14\mid 23)$,
which we call the {\em special boundary divisors}.
(We will also use this term to refer to individual components of $D(12\mid
34)$, etc.)
\par
To obtain the Gromov-Witten invariants defined in {\S}\ref{GWinv},
we will not need to analyze these divisors completely; it will be
enough to identify just those components that affect our calculations.
A divisor $D$ on $\lmodnd$
is called \emph{numerically irrelevant} if
\begin{equation}\label{numirr}
\int e_1^{*}(\gamma_1) \cup \cdots \cup e_n^{*}(\gamma_n) \cap [D] = 0
\end{equation}
for all $\gamma_1,\ldots,\gamma_n \in A^*(S_2)$.  Divisors that are
not numerically irrelevant are said to be \emph{numerically relevant}.
Similarly, for $n \geq 4$, a divisor is called
\emph{irrelevant with respect to the base} if
$$
\int e_1^{*}(\gamma_1) \cup \cdots \cup e_n^{*}(\gamma_n) \cap [D] = 0
$$
for all $\gamma_1,\gamma_2,\gamma_3,\gamma_4 \in A^*(\pp)$ and all
$\gamma_5,\ldots,\gamma_n \in A^*(S_2)$.
\par
Let
$\tau\colon  \lmodnd \to \ppmod{0}{d}$ be the morphism that composes a
stable map $\mu\colon  C \to S_2$ with the projection $S_2 \to \pp$, forgets
all markings, and contracts any twigs
that have become unstable.
In analogy with Proposition 5.4 of
\cite{ErnstromKennedy}, we have the following result.

\begin{prop}\label{irrellemma}
A divisor $D$ on the stack $\lmodnd$ is numerically irrelevant if
$\cod \tau(D) \geq 2$.
\end{prop}

\begin{proof}
Without loss of generality, we may assume that $D$ is irreducible.
By linearity we may assume that the $\gamma_j$ in (\ref{numirr}) come
from the following basis for $A^{*}(S_2)$:
$$
\{ 1,h,\hd,h^2,\hd^2, h^2 \hd, i,hi,\hd i, h^2i, \hd^2i, h^2\hd i \}.
$$
Note first that if any $\gamma_j$, say $\gamma_n$, is the identity
element then, by the projection formula applied to the forgetful
morphism $\pi\colon  \lmodnd \to \lmod{n-1}{d}$, we
have
$$
\int e_1^{*}(\gamma_1) \cup \cdots \cup e_{n-1}^{*}(\gamma_{n-1})
   \cup e_n^{*}(1) \cap [D]
=  \int e_1^{*}(\gamma_1) \cup \cdots \cup e_{n-1}^{*}(\gamma_{n-1})
   \cap \pi_{*}[D],
$$
since the evaluation maps $e_1,\ldots,e_{n-1}$ factor through $\pi$.
But $\pi_*[D]$ is zero unless $\dim \pi(D) = \dim D = \dim \lmod{n-1}{d}$.
By the irreducibility of the moduli stack, we
therefore have $\pi(D) = \lmod{n-1}{d}$ and hence that
$\tau(D) = \ppmod{0}{d}$.  Given that $\cod \tau(D)$ is at least
$2$, we see that the proposition holds in this instance.

\par
Next, consider, for $k \leq n$, the class
$$
\alpha_k = e_1^{*}(\gamma_1) \cup \cdots \cup e_k^{*}(\gamma_k) \cap [D]
$$
where none of the $\gamma_j$'s is the identity.  We claim that if
$\alpha_k$ is nonzero, then it can be represented by a cycle $D_k$ for
which
\begin{equation}\label{codeqn}
\cod \tau(D_k) \geq 2 + \sum_{j=1}^{k}\left( \cod \gamma_j - 1\right).
\end{equation}

\par
To establish the claim, we use induction on $k$, the case $k = 0$
being trivial.  In view of the remarks above, we
need only consider the cases where $\cod \gamma_k \geq 2$.  For
$\gamma_k = h^2, \hd^2$, or $h^2\hd$, the argument is
given in the proof of Proposition 5.4 of \cite{ErnstromKennedy}.
Thus we prove the result when
$\gamma_k = hi$, $h^2i$, $\hd i$, $\hd^2i$, or $h^2\hd i$.

\begin{enumerate}
\item
Suppose $\gamma_k = hi$.  If $L \subset \pp$ is a general line,
then $\alpha_k$ may be represented by a cycle $D_k$ for which
$\tau(D_k)$ is contained in the intersection of $\tau(D_{k-1})$ and
the set of curves with a singularity on $L$.  Hence
$\cod \tau(D_k) > \cod \tau(D_{k-1})$ and (\ref{codeqn}) holds by induction.

\item
Suppose $\gamma_k = h^2i$.  If $P \in \pp$ is a general point,
then $\alpha_k$ may be represented by a cycle $D_k$ for which
$\tau(D_k)$ is contained in the intersection of $\tau(D_{k-1})$ and
the set of curves with a singularity at $P$.  Hence
$\cod \tau(D_k) \geq 2+\cod \tau(D_{k-1})$ and (\ref{codeqn}) holds by
induction.

\item
Suppose $\gamma_k = \hd i$.  Let $P \in \pp$ be a general point.  Since
$i = z + 3(h-\hd)$, we have that $\alpha_k$ is represented by a cycle
$D_k$ for which $\tau(D_k)$ is in the intersection of $\tau(D_{k-1})$
and the union of the set of curves having a flex tangent that passes
through $P$ and the set of curves whose tangent cone at a singular
point passes through $P$.  Thus $\cod \tau(D_k) > \cod \tau(D_{k-1})$.

\item
Suppose $\gamma_k = \hd^2i$.  Let $L \subset \pp$ be a general line.
Using $z$ in place of $i$ as in case 3, we have that $\alpha_k$ is
represented by a cycle $D_k$ for which $\tau(D_k)$ is in the
intersection of $\tau(D_{k-1})$ and the union of the set of curves
having $L$ as flex tangent and the set of curves whose tangent cone at
a singular point contains $L$.  Thus $\cod \tau(D_k) \geq 2+ \cod
\tau(D_{k-1})$.

\item
Suppose $\gamma_k = h^2\hd i$.  Let $(P,L)$ be a general flag.  Using
$z$ in place of $i$, we have that $\alpha_k$ is represented by a cycle
$D_k$ for which $\tau(D_k)$ is in the intersection of $\tau(D_{k-1})$
and the union of the set of curves with a flex at $P$ and flex
tangent equal to $L$ and the set of curves singular at $P$ whose tangent cone
contains $L$.  Thus $\cod \tau(D_k) \geq 3+ \cod \tau(D_{k-1})$.

\end{enumerate}

\par
To finish the proof, note that
$\int e_1^{*}(\gamma_1) \cup \cdots \cup e_n^{*}(\gamma_n) \cap [D]$
must be zero unless $\sum_{j=1}^{n} \cod \gamma_j = 3d-2+n$, i.e.,
unless $2+\sum_{j=1}^{n} (\cod \gamma_j-1) = 3d$.
But in this case (\ref{codeqn}) implies that
$$
\cod \tau(D_n) \geq 3d > \dim \ppmod{0}{d}.
$$
Hence $D_n$ must be zero, as desired.
\end{proof}

\par

\begin{cor}\label{irrelcor}
If $D$ is a numerically relevant special boundary divisor on the stack
$\lmodnd$, then $\tau(D)$ is one of
the following:
\begin{itemize}
   \item[(1)]  all of $\ppmod{0}{d}$,
   \item[(2)]  a divisor $D(\emptyset,\emptyset;d_1,d_2)$, where
   $d_1$, $d_2$ are positive and $d_1+d_2=d$,
   \item[(3)]  a divisor $\cuspmod{d}$, for $d \geq 2$.
\end{itemize}
\end{cor}

\begin{proof}
An immersion $\pl \to \pp$ has a unique (unmarked) second-order stable
lift, whose
source is again $\pl$.
Thus the general $n$-marked stable lift of an
immersion likewise has $\pl$ as source. Thus the image in
$\ppmod{0}{d}$ of a numerically relevant divisor is
either the entire stack or a divisor representing
non-immersive maps. Every non-immersion either has a reducible source
or is singular at one or more points; thus it is represented by a point of
one of the listed divisors.
\end{proof}

\par
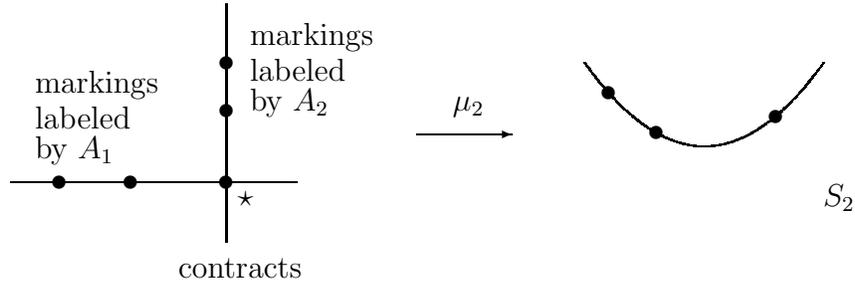
\begin{figure}
   \begin{center}
   \setlength{\unitlength}{0.125in}
      \begin{picture}(36,12)
        %left part
        \thicklines
        \put(9,1.5){\line(0,1){10}}
        \thinlines
        \put(0,4){\line(1,0){12}}
        \put(2,4){\circle*{.5}}
        \put(5,4){\circle*{.5}}
        \put(9,4){\circle*{.5}}
        \put(9,7){\circle*{.5}}
        \put(9,9){\circle*{.5}}
        \put(9.5,3){$\star$}
        \put(7,0){contracts}
        \put(1,5){\shortstack[l]{markings\\ labeled\\ by $A_1$}}
        \put(10,7){\shortstack[l]{markings\\ labeled\\ by $A_2$}}
        \put(17,6){\vector(1,0){4}}
        \put(18.5,7){$\mu_2$}
        %right part
        \qbezier(24,9)(29,2)(34,9)
        \put(25,7.75){\circle*{.5}}
        \put(27,6.1){\circle*{.5}}
        \put(32,6.75){\circle*{.5}}
        \put(34,3){$S_2$}
      \end{picture}
   \end{center}
   \caption{A second-order stable lift over an (unmarked) immersion.}
   \label{collapsepic}
   \end{figure}

Suppose that $\mu$ is an element of $\lmodnd$ for which $\tau(\mu)$
is an immersion. According to Proposition
\ref{fibersquare}, one of the twigs is the strict lift of
$\tau(\mu)$. Other twigs, if any, must be mapped to single points of
$S_2$, and thus must carry at least two markings.
Such a stable lift is
represented by a point of $D(12\mid 34)$ if and only if the markings
$1$ and $2$ are separated from the markings $3$ and $4$ by at least
one attachment point. Thus there is a component
$$
\lmod{A_1 \cup \{\star \}}{d} \times_{S_2} \sempmod{A_2 \cup \{\star \}}{0}
$$
of $D(12\mid 34)$
corresponding to every partition in which $1$, $2$ belong to one
of the subsets and $3$, $4$ to the other. Here $\star$ indicates the
marking used to create the fiber product; $A_1$ indexes the markings
on the strict lift twig and $A_2$ those on the constant twig. (A
typical member is shown in Figure \ref{collapsepic}.)
%Figure for the stable lift for one of the relevant divisors

\par
In Proposition \ref{node}, we have described the unique (unmarked)
second-order stable lift of a general element of
$D(\emptyset,\emptyset;d_1,d_2)$.
(Also see Figure \ref{nodepic}.)
According to Proposition
\ref{fibersquare}, on an $n$-marked stable lift, the markings may be
distributed on this five-twig curve in any way (with additional twigs
arising as markings are brought together).  Thus there is a component
$D^2(A_1,A_2,A_3,A_4,A_5;d_1,d_2)$ of $D(12\mid 34)$ corresponding to
each valid partition, i.e., in which marks $1$ and $2$ occur are separated
from marks $3$ and $4$ by at least one attachment point. Here $A_1$
indexes the markings on the peripheral twig carrying the strict lift
of the map of degree~$d_1$, and
$A_2$ indexes the markings on the adjacent twig, etc. To avoid
redundancy, we may assume that marks $1$ and $2$ occur before marks
$3$ and $4$. The divisor $D^2(A_1,A_2,A_3,A_4,A_5;d_1,d_2)$ may be
expressed as a fiber product
\begin{equation*}
\begin{split}
\lmod{A_1 \cup \{\clubsuit\}}{\,&d_1}
\times_{S_2}
\ramstack{A_2 \cup \{\clubsuit,\diamondsuit\}}{2,3}
\times_{S_2}
\ramstack{A_3 \cup \{\diamondsuit,\heartsuit\}}{1,2} \\
&\times_{S_2}
\ramstack{A_4 \cup \{\heartsuit,\spadesuit\}}{2,3}
\times_{S_2}
\lmod{A_5 \cup \{\spadesuit\}}{d_2}
\end{split}
\end{equation*}
in which the first and last factors are stacks of stable lifts.  The
middle factor is the stack whose general member represents a double
cover of the lift to $S_2$ of a fiber of $S_1$ over a point of $\pp$,
ramified at the two markings $\diamondsuit$ and $\heartsuit$ and
carrying additional markings indexed by $A_3$.  Similarly, the second
and fourth factors are stacks whose general member represents a triple
cover of a fiber of $S_2$ over a point of $S_1$, ramified only at the
two special markings.  The fiber products are created using the
evaluation maps at $\clubsuit$, $\diamondsuit$, $\heartsuit$, and
$\spadesuit$.
\par
Similarly,
in Proposition \ref{cusp}, for $d \geq 3$, we have described the unique
unmarked
second-order lift of a general element of $\cuspmod{d}$,
and on an
$n$-marked stable lift the markings may be distributed in any way.
(Also see Figure \ref{cusppic}.)
Thus there is a component $C^2(A_1,A_2,A_3,A_4,A_5;d)$ of $D(12\mid 34)$
corresponding to each valid partition.  Likewise, there is a
component $C^2(A_1,A_2,A_3,A_4,A_5,A_6,A_7;2)$ for each valid
partition.
\par
In view of these remarks and Corollary \ref{irrelcor}, the following
theorem is now clear.

\begin{thm}\label{relevant}
If $D$ is a numerically relevant component of $D(12\mid 34)$, then it
is one of the following:
\begin{itemize}
  \item[(1)] a divisor $\lmod{A_1 \cup \{\star \}}{d} \times_{S_2}
               \sempmod{A_2 \cup \{\star \}}{0}$ where $|A_2| \geq 2$;
  \item[(2)] a divisor $D^2(A_1,A_2,A_3,A_4,A_5;d_1,d_2)$;
  \item[(3)] a divisor $C^2(A_1,A_2,A_3,A_4,A_5;d)$ with $d \geq 3$;
  \item[(4)] a divisor $C^2(A_1,A_2,A_3,A_4,A_5,A_6,A_7;2)$.
\end{itemize}
\end{thm}

\par
\begin{thm}\label{baseirrel}
The numerically relevant components of type (3) or (4)
are irrelevant with respect to the base. Among those of type (2),
only those in which markings 1 and 2 belong to $A_1$, and markings 3
and 4 belong to $A_5$, are relevant
with respect to the base.
\end{thm}
\begin{proof}
Consider a component $D$ of type (2) in which the the first marking
belongs to $A_2$. Then the evaluation map $e_1$
factors through projection of the fiber product
$D$ onto its second factor, the stack
$\ramstack{A_2 \cup \{\clubsuit,\diamondsuit\}}{2,3}$.
Since this stack represents maps to fibers of $S_2$ over $S_1$,
there is a morphism which picks out the image of the fiber and thus a
composite morphism $\epsilon\colon  D \to \pp$.
If $\gamma_1$ is an element of
$A^*(\pp)$, the classes $e_1^{*}(\gamma_1)$ and $\epsilon^{*}(\gamma_1)$
agree.
Let $D'$ denote $D^2(A_1,A_2 \setminus \{1\},A_3,A_4,A_5;d_1,d_2)$, and
let $\epsilon'\colon  D' \to \pp$ be the analogous morphism. Then $\epsilon =
\epsilon'
\circ \varphi$, where $\varphi\colon D \to D'$ is the forgetful morphism.
If $\gamma_2,\dots,\gamma_n$ are elements of $A^*(S_2)$,
then by the projection formula
\begin{equation*}
\begin{split}
\int e_1^{*}(\gamma_1) \cup e_2^{*}(\gamma_2) \cup \cdots \cup
e_n^{*}(\gamma_n) \cap [D]
& =\int \epsilon^{*}(\gamma_1) \cup e_2^{*}(\gamma_2) \cup \cdots \cup
e_n^{*}(\gamma_n) \cap [D] \\
& =\int {\epsilon'}^{*}(\gamma_1) \cup e_2^{*}(\gamma_2) \cup \cdots \cup
e_n^{*}(\gamma_n) \cap \varphi_*[D].
\end{split}
\end{equation*}
Since the relative dimension of $D$ over $D'$ is 1, this class vanishes.
Hence $D$ is irrelevant with respect to the base.
\par
A similar argument applies if the second, third, or fourth marking
belongs to $A_2$. The argument also works if any of these markings
belong to $A_3$ or $A_4$, since again the corresponding evaluation
map factors through projection onto a factor of $D$ representing maps
to fibers of $S_2$ over $S_1$, or representing maps to fibers of $S_1$ over
$\pp$. For a component of type (3) or (4) we may also employ the same sort
of argument, since for any such component at least one (in fact two) of
the four special markings must lie on a twig which either maps to a fiber
of $S_2$ over $S_1$, or maps to a fiber of $S_1$ over $\pp$, or is
contracted to a point.
\end{proof}
\par
Having identified the components of $D(12\mid 34)$ which are relevant
with respect to the base, we must investigate their multiplicities.
Here it is important to distinguish between the stack $\lmodnd$ and
its coarse moduli space.   A general point of the closed substack
$D^2(A_1,\emptyset,\emptyset,\emptyset,A_5;d_1,d_2)$ represents a stable
map with an
automorphism group of order $18$;
thus the map from this substack to
its image on the coarse moduli space---a divisor---has degree~$1/18$.
Under the birational morphism from the coarse moduli space of $\lmodnd$
to that of $\ppmod{n}{d}$, this divisor is mapped onto a divisor. Now
on the coarse moduli space of $\ppmod{n}{d}$, all components of $D(12\mid
34)$ appear with multiplicity $1$. Hence on the stack $\lmodnd$ the
fundamental
class of the substack
$D^2(A_1,\emptyset,\emptyset,\emptyset,A_5;d_1,d_2)$ appears in $D(12\mid
34)$ with multiplicity $18$.
By Proposition \ref{fibersquare}, a similar statement holds when
there are markings on the middle twigs:
$D^2(A_1,A_2,A_3,A_4,A_5;d_1,d_2)$ appears in $D(12\mid 34)$ with multiplicity
$18$. The same sort of argument shows that a divisor of type (1)
(the general point of which represents an automorphism-free stable
map) occurs in $D(12\mid 34)$ with multiplicity $1$.
(See \cite{Vistoli} for a general discussion of intersection theory on
stacks.)
\par

\section{Potentials and partial differential equations}
\label{potpde}

\par
We now present potential functions associated to each of the special
boundary divisors that are relevant with respect to the base;
these are the potentials that will be needed in order to recursively
calculate the second-order Gromov-Witten invariants.
Since many of the intermediate expressions are exceedingly long,
we will describe them rather than write them out explicitly.
In these calculations we employ two ordered bases for
$A^*(S_2)$, which are dual under the intersection pairing.
Here is the ``$z$-basis'':

\begin{center}
{\renewcommand{\arraystretch}{1.5}
\begin{tabular}{c|c|c|c|c|c|c|c|c|c|c|c}
$Y_{000}$ & $Y_{100}$ & $Y_{200}$ & $Y_{010}$ & $Y_{020}$ & $Y_{210}$ &
$Y_{001}$ & $Y_{101}$ & $Y_{201}$ & $Y_{011}$ & $Y_{021}$ & $Y_{211}$ \\

\hline
$1$ & $h$ & $h^2$ & $\hd$ & $\hd^2$ & $h^2\hd$ &
$z$ & $hz$ & $h^2 z$ & $\hd z$ & $\hd^2 z$ & $h^2\hd z$
\end{tabular}}
\end{center}

\noindent
This is the ``$i$-basis'':

\begin{center}
{\renewcommand{\arraystretch}{1.5}
\begin{tabular}{c|c|c|c|c|c|c|c|c|c|c|c}
$Y_{211}^*$ & $Y_{021}^*$ & $Y_{011}^*$ & $Y_{201}^*$ & $Y_{101}^*$ &
$Y_{001}^*$ &
$Y_{210}^*$ & $Y_{020}^*$ & $Y_{010}^*$ & $Y_{200}^*$ & $Y_{100}^*$ &
$Y_{000}^*$ \\ \hline
$h^2\hd i$ & $\hd^2 i$ & $\hd i$ & $h^2 i$ & $hi$ & $i$ &
$h^2 \hd$ & $\hd^2$ & $\hd$ & $h^2$ & $h$ & $1$
\end{tabular}}
\end{center}

\noindent
Here $z$ is the class of $Z$ and $i$ is
the class of the divisor at infinity $I$.  Note that the subscript
$k$ on $Y_k$ is a vector having three integer
components, and that the basis element dual to $Y_k$ is $Y_{k^*}^{*}$,
where $k+k^* = 211$ or $121$.  We will write a general element of
$A^*(S_2)$ as $\sum y_k Y_k$ or as $\sum y_l^{*} Y_l^{*}$.

\par
We write the classical potential $\calP$ of $S_2$ in a somewhat
nonstandard way by using both bases.
Let $\gamma = \sum y_k Y_k$ and $\delta = \sum y_l^{*} Y_l^{*}$. Then
\begin{align*}
  \calP &= \sum_{n \geq 2} \frac{1}{n!} \int e_1^{*}(\gamma)
           \cup \cdots \cup e_n^{*}(\gamma) \cup e_{n+1}^{*}(\delta)
           \cap [\sempmod{n+1}{0}]\\
        &= \frac{1}{2} \int \gamma \cup \gamma \cup \delta \cap [S_2].
\end{align*}
\par
The quantum potential of second-order stable lifts is defined by
$$
\calN = \sum_{\substack{d \geq 1\\ n \geq 0}}\frac{1}{n!}
         \int e_1^{*}(\gamma) \cup \cdots \cup e_n^{*}(\gamma)
              \cap [\lmodnd].
$$
Using the basic properties of Gromov-Witten invariants given in
{\S}\ref{GWinv}, we may express $\calN$ as
$$
\calN = \sum_{d \geq 1}N^{(d)}E^{(d)}
$$
where $E^{(d)} = \exp(dy_{100}+(2d-2)y_{010}+(3d-6)y_{001})$ and
$N^{(d)}$ is a polynomial in the variables $y_k$ corresponding to
the eight elements of the $z$-basis other than the identity and
the divisors. The polynomial $N^{(d)}$ is weighted
homogeneous, where $\operatorname{weight}(y_k) = \cod(Y_k) -1$; it is
the generating series of Gromov-Witten invariants, without divisor
conditions, of degree~$d$.  To be more explicit, let
$$
\mathbf{a} =
(a_{200},a_{020},a_{210},a_{101},a_{201},a_{011},a_{021},a_{211})
$$
be a vector of nonnegative integers.  If $k$ is any of the possible
subscripts appearing in $\mathbf{a}$, let $|k|$ denote the sum of its
entries, and set
$$
\|\mathbf{a}\| = \sum_k (|k|-1)a_k.
$$
Let $\mathbf{y}^{\mathbf{a}} = \prod_{k} y_k^{a_k}$ and
$\mathbf{a}! = \prod_{k} a_{k}!$.  Set
$N_d(\mathbf{a}) = \langle \prod_i \gamma_i\rangle_d$ where $a_{200}$
of the $\gamma_i$ are $Y_{200}$, $a_{020}$ are $Y_{020}$, etc.  Then
$$
N^{(d)} = \sum_{\|\mathbf{a}\| = 3d-1}
          N_d(\mathbf{a}) \frac{\mathbf{y}^{\mathbf{a}}}{\mathbf{a}!}.
$$

\par
Next we introduce a potential $\calT$ for the maps on the three middle
twigs of the divisor $D^2(A_1,A_2,A_3,A_4,A_5;d_1,d_2)$:
$$
\calT = \sum_{n \geq 0} \frac{1}{2n!}\int e_1^{*}(\gamma) \cup \cdots
\cup
e_n^{*}(\gamma) \cup e_{n+1}^{*}(\delta) \cup e_{n+2}^{*}(\delta)
\cap [\tail_n]
$$
Here $\tail_n$ is the stack whose general point represents a
stable map into $S_2$ with a three-twig source:  the central twig maps as
a double cover of the lift of a fiber of $S_1$ over $\pp$ with
ramification occurring at the attachment points, and the peripheral
twigs each map as a triple cover of a fiber of $S_2$ over $S_1$ with
ramification concentrated at two points, one of which is an attachment
point.  In the definition above, we express both $\gamma$ and
$\delta$ in terms of the $z$-basis, and we use the fundamental class
of the stack rather than the fundamental class of its
associated coarse moduli space. (As explained at the end of
{\S}\ref{sbd}, the two classes differ by a factor of $18$.)
\par
The stack $\tail_n$ is the fiber product of three simpler stacks, each
with an associated potential, which we now describe. Let
$\ramstack{n,2}{1,2}$ denote the stack of stable maps whose general
member represents an $(n+2)$-marked map $\pl \to S_2$ satisfying
these conditions:
\begin{itemize}
 \item the map is a double cover of a fiber of the lift to $S_2$ of a
 fiber of $S_1$ over $\pp$;
 \item the last two markings always occur at the ramification points.
\end{itemize}
In other words, this stack represents central twigs of members of
$\tail_n$, with two extra markings for ``gluing.''
Again we emphasize that we are employing the fundamental class of the
stack rather than its associated coarse moduli space. (Here the two
classes differ by a factor of $2$.)
Let $\calR(1,2)$ denote the associated potential:
$$
\calR(1,2) = \sum_{n \geq 0} \frac{1}{2n!}\int e_1^{*}(\gamma) \cup
\cdots \cup
e_n^{*}(\gamma) \cup e_{n+1}^{*}(\delta) \cup e_{n+2}^{*}(\delta)
\cap [\ramstack{n,2}{1,2}].
$$
We will express both $\gamma = \sum y_{k}Y_{k}$ and $\delta = \sum
z_{k}Y_{k}$
in terms of the $z$-basis.
\par
Similarly, we define $\ramstack{n,2}{2,3}$ to be the stack for the
peripheral twigs. The general member represents an $(n+2)$-marked map
$\pl \to S_2$ that is a triple cover of a fiber of $S_2$ over $S_1$,
and again we demand that the last two markings occur at the ramification
points.  The associated potential is
$$
\calR(2,3) = \sum_{n \geq 0} \frac{1}{2n!}\int e_1^{*}(\gamma) \cup
\cdots \cup
e_n^{*}(\gamma) \cup e_{n+1}^{*}(\delta) \cup e_{n+2}^{*}(\delta)
\cap [\ramstack{n,2}{2,3}].
$$
In $\calR(2,3)$ we will express $\gamma = \sum y_{k}Y_{k}$ in terms of
the $z$-basis and $\delta = \sum z_{l}^{*}Y_{l}^{*}$ in terms of the
$i$-basis.
\par
\begin{lem}\label{r12lem}
Let $y_{110}=y_{020}$ and let $z_{110}=z_{020}$.
The potential $\calR(1,2)$ is the sum of terms in
$$
\frac{1}{2}\exp(2y_{010})\exp(2y_{110})\exp(2y_{210})\exp(z_{010})
\exp(z_{110})\exp(z_{210})
$$
that are quadratic in $z_{010},z_{110},z_{210}$ and for which the
(vector) sum of subscripts is $(2,n+2,0)$ for some $n$.
\end{lem}

\begin{proof} First use the alternative basis to the $z$-basis where $h\hd$
replaces $\hd^2$ and $h\hd z$
replaces $\hd^2 z$. Denote the coefficients of $\calR(1,2)$ with respect to
this
alternative basis by $\tilde{y}_k$'s and $\tilde{z}_l$'s.  Each basis
element may be written in the form $h^a\gamma$, where $a = 0$, $1$,
or $2$, and $\gamma$ does not involve $h$.  The Gromov-Witten
invariant
\begin{equation*}
\begin{split}
\int e_1^{*}(h^{a_1}\gamma_1) \cup \cdots
     &\cup e_n^{*}(h^{a_n}\gamma_n) \\
     &\cup e_{n+1}^{*}(h^{a_{n+1}}\gamma_{n+1}) \cup
     e_{n+2}^{*}(h^{a_{n+2}}\gamma_{n+2}) \cap [\ramstack{n,2}{1,2}]
\end{split}
\end{equation*}
is zero unless $a_1 + \cdots + a_{n+2} = 2$, and otherwise equals
$$
\int e_1^{*}(\gamma_1) \cup \cdots \cup e_n^{*}(\gamma_n)
     \cup e_{n+1}^{*}(\gamma_{n+1}) \cup
     e_{n+2}^{*}(\gamma_{n+2}) \cap [L]
$$
where $[L]$ denotes the fundamental class of the lift of a fiber of
$S_1$
over $\pp$.  This in turn must be zero if $z$ appears in any
$\gamma_i$ (since the divisor at infinity $I$ is disjoint from the
dual divisor $Z$).  It also vanishes unless all $\gamma_i$'s equal
$\hd$, in which case its value is $2^{n-1}$.  Thus
$$
\calR(1,2) = \sum_{n, \mathbf{k}}\frac{1}{2n!} \,C\tilde{y}_{k_1}
\cdots \tilde{y}_{k_n}\tilde{z}_{k_{n+1}}\tilde{z}_{k_{n+2}}
$$
where the coefficient $C$ is $2^{n-1}$ if the first entry of $k_1+\cdots
+k_{n+2}$ is $2$, the last entry of this sum is $0$, and each second
entry is $1$; otherwise the coefficient is zero.  Thus $\calR(1,2)$
is the sum of terms in
$$
\frac{1}{2}\exp(2\tilde{y}_{010})\exp(2\tilde{y}_{110})\exp(2\tilde{y}_{210})
\exp(\tilde{z}_{010})\exp(\tilde{z}_{110})\exp(\tilde{z}_{210})
$$
that are quadratic in $\tilde{z}_{010},\tilde{z}_{110},\tilde{z}_{210}$
and for which the sum of subscripts is $(2,n+2,0)$ for some $n$.  To
obtain the lemma,
note that $y_k = \tilde{y}_k$ for all subscripts $k$ appearing in the
formula (which do not include the subscripts $200$ or $201$).
\end{proof}

By a similar argument, we establish the following lemma.  Recall that
we use both the $z$-basis and the $i$-basis to write the
potential $\calR(2,3)$.

\begin{lem}\label{r23lem}
The potential $\calR(2,3)$ is the sum of terms in
$$
\frac{1}{3}\prod_{k_3 =1}\exp(3y_k)\prod_{l_3 =1}\exp(y_l^{*})
$$
that are quadratic in the $y_l^{*}$'s and for which the
(vector) sum of subscripts is $(2,1,n+2) $or $(1,2,n+2)$ for some $n$.
\end{lem}

\begin{proof}
The Gromov-Witten invariant
\begin{equation*}
\begin{split}
\int e_1^{*}(h^{a_1}\hd^{b_1}\gamma_1) \cup \cdots
          &\cup e_n^{*}(h^{a_n}\hd^{b_n}\gamma_n) \\
          &\cup e_{n+1}^{*}(h^{a_{n+1}}\hd^{b_{n+1}}\delta_{n+1}) \cup
     e_{n+2}^{*}(h^{a_{n+2}}\hd^{b_{n+2}}\delta_{n+2}) \cap
     [\ramstack{n,2}{2,3}]
\end{split}
\end{equation*}
is zero unless $a_1 + \cdots + a_{n+2} = 1 \text{ or } 2$, $b_1 + \cdots +
b_{n+2} =
2 \text{ or } 1$ (respectively), $\gamma_1 = \ldots = \gamma_n = z$,
and $\delta_{n+1}=\delta_{n+2} = i$, in which case it has value
$3^{n-1}$.  Thus
$$
\calR(2,3) = \sum_{n, \mathbf{k}}\frac{1}{2n!} \,C y_{k_1}
\cdots y_{k_n}y_{k_{n+1}}^{*}y_{k_{n+2}}^{*}
$$
where the coefficient $C$ is $3^{n-1}$ if the first entry of $k_1+\cdots
+k_{n+2}$ is $1 \text{ or } 2$, the second entry is $2 \text{ or } 1$
(respectively), and the last entry of each
$k_j$ is $1$; it is zero otherwise.  Hence the lemma follows.
\end{proof}
\par
To put the potentials $\calR(1,2)$ and $\calR(2,3)$ together to obtain
an expression for the potential $\calT$, we appeal to some
general results. We will need the following easy consequence of the
chain rule.

\begin{prop}\label{psprop}
Suppose $\{ Y_0,\ldots,Y_m \}$ and $\{ Z_0,\ldots,Z_m \}$ are
two ordered bases for a vector space.
Let $\gamma = \sum_{k=0}^{m}y_k Y_k$ and $\delta = \sum_{l=0}^{m}z_l
Z_l$.  Suppose $\calK$ is a power series in
$y_0,\ldots,y_m,z_0,\ldots,z_m$ that can be
written as
$$
\calK = \sum_{a,b \geq 0} \frac{K_{a,b}(\gamma^{a};\delta^{b})}{a!b!}
$$
where each $K_{a,b}$ is multilinear and symmetric in both the first
$a$ and last $b$ arguments.  Let $k_1,\ldots, k_r$ and
$l_1,\ldots, l_s$ have values in $\{0,\ldots,m\}$.  Then
$$
\frac{\partial^{r+s}\calK}
{\partial y_{k_1}\cdots \partial y_{k_r}\partial z_{l_1}\cdots \partial
z_{l_s}}
= \sum_{a,b \geq 0}\frac{K_{a+r,b+s}(\gamma^{a}Y_{k_1} \cdots
Y_{k_r};\delta^{b}Z_{l_1} \cdots Z_{l_s})}{a!b!}.
$$
\end{prop}

\par
Now suppose that $X$ is an arbitrary scheme, that
$\beta_1+\beta_2$ is a partition of the element $\beta \in A_{1}X$, and that
$A_1 \cup A_2$ is a partition of the index set $\{1,\dots,n\}$.
Suppose that $M_1$ is a substack of
${\overline M}_{A_1\cup\{\star\}}(X,\beta_1)$
and $M_2$ is a substack of
${\overline M}_{A_2\cup\{\star\}}(X,\beta_2)$.
Then the fiber product
$M_1 \times_{X} M_2$, formed using the evaluation maps at
$\star$, is a substack of ${\overline M}_n(X,\beta)$.
We have a fiber square
$$
\begin{CD}
M_1 \times_{X} M_2 @>\iota>> M_1 \times M_2 \\
@VVV        @VVV\\
X^{n+1} @>D>> X^{n+2}
\end{CD}
$$
where $\iota$ is inclusion, $D$ is the diagonal inclusion that
repeats the last
coordinate, and the coordinates of the vertical maps are evaluation maps to
$X$.
In particular the evaluation maps $e_{n+1},e_{n+2}$ on the right are at the
two
markings labeled by $\star$; note
that they agree when restricted to the fiber product.
\par
Let $\{ Y_0,\ldots,Y_m \}$ and $\{ Y_m^{*},\ldots,Y_0^{*} \}$ be
ordered bases for $A^*(X)$ that are dual with respect to the
intersection pairing.
Then the fundamental class of the diagonal
$\Delta$ in $X \times X$ is
$$
[\Delta] = \sum_{k=0}^{m}Y_k \times Y_{m-k}^{*}.
$$
Thus
\begin{equation*}
  \iota_*\left( e_1^{*}(\gamma_1) \cup \cdots \cup e_n^{*}(\gamma_n)
\right)
    = \sum_{k=0}^{m}e_1^{*}(\gamma_1) \cup \cdots \cup
    e_n^{*}(\gamma_n) \cup e_{n+1}^{*}(Y_k) \cup
    e_{n+2}^{*}(Y_{m-k}^{*}).
\end{equation*}
Taking degrees, we obtain the following proposition.

\begin{prop}\label{fibprodprop}
In this situation
\begin{equation*}
\begin{split}
\int & e_1^{*}(\gamma_1) \cup \cdots \cup e_n^{*}(\gamma_n) \cap[M_1
\times_{X} M_2] =\\
      & \sum_{k=0}^{m}\left(
      \int \bigcup_{t \in A_1}e_t^{*}(\gamma_t) \cup e_{n+1}^{*}(Y_k)
      \cap [M_1]\right)
      \left(
      \int \bigcup_{t \in A_2}e_t^{*}(\gamma_t) \cup
e_{n+2}^{*}(Y_{m-k}^{*})
      \cap [M_2]\right).
\end{split}
\end{equation*}
\end{prop}

\par
Let $\gamma = \sum_{k=0}^{m}y_k Y_k$.
Suppose, for $i = 1,2$ that the stack $M_i$ has potential $\calM_i$, i.e.,
$$
\calM_i =
\sum_{n \geq 0} \frac{1}{n!} \int e_1^{*}(\gamma) \cup \cdots
                                 \cup e_n^{*}(\gamma) \cap[M_i].
$$
\begin{prop}\label{prodpot}
The potential $\calM$ for the fiber product $M_1 \times_{X} M_2$ is
$$
\calM =
\sum_{k=0}^{m} \frac{\partial \calM_1}{\partial y_k}
               \frac{\partial \calM_2}{\partial y_{m-k}^{*}}.
$$
\end{prop}

\begin{proof}
Set $\gamma_1 = \ldots = \gamma_n = \gamma  = \sum y_k Y_k$.  Apply
Proposition \ref{fibprodprop} and sum over all $n \geq 0$.  Then apply
Proposition \ref{psprop}.
\end{proof}

\par
Using Proposition
\ref{prodpot} twice, we deduce the following result.

\begin{prop}\label{tailpot}
The potential $\calT$ is given by
$$
\calT = \frac{1}{2}\sum_{s,t}\frac{\partial \calR(2,3)}{\partial
y_{s^*}^{*}}
                             \frac{\partial^2 \calR(1,2)}{\partial y_s
\partial y_t}
                             \frac{\partial \calR(2,3)}{\partial
y_{t^*}^{*}}.
$$
\end{prop}
\noindent
Proposition \ref{tailpot}, together with Lemmas \ref{r12lem} and
\ref{r23lem}, is sufficient to give an explicit expression for
$\calT$.
\par
Next we derive partial differential equations in a manner analogous to the
derivation of the WDVV-equation
of quantum cohomology (equation (58) of \cite{FultonP}). Fix an 
ordered quadruple $Y_i, Y_j, Y_k, Y_l$
of elements from the $z$-basis, and set
$$
\calG(ij \mid kl) = \sum_{\substack{n \geq 4\\ d \geq 0}}
\frac{1}{n!}
  \int e_1^{*}(Y_i) \cup e_2^{*}(Y_j) \cup e_3^{*}(Y_k) \cup
e_4^{*}(Y_l)
      \cup e_5^{*}(\gamma) \cup \cdots \cup e_n^{*}(\gamma)
      \cap [D(1 2 \mid 3 4)],
$$
where $D(1 2 \mid 3 4)$ denotes the special boundary divisor on
$\ppmod{n}{d}$. Note that the same notation
$D(1 2 \mid 3 4)$ is used for divisors in the various stacks
$\ppmod{n}{d}$.

\par
The linear equivalence $D(1 2 \mid 3 4) \sim D(1 4 \mid 2 3)$
implies the following equality.
\par
\begin{thm}\label{Gijklprop} For all subscripts $i$, $j$, $k$, and $l$,
$$
\calG(i j \mid k l) = \calG(i l \mid j k).
$$
\end{thm}
\noindent
In general this will be a complicated identity, since it will
involve contributions from all numerically relevant components
of $D(1 2 \mid 3 4)$ and $D(1 4 \mid 2 3)$.
However, if
each subscript is either $100$ or $200$ (i.e., if the classes
$Y_i$, $Y_j$, $Y_k$, $Y_l$ are pullbacks to $S_2$ of classes from
$\pp$), then we need to know only the components which are relevant
with respect to the base. In these cases,
Theorem \ref{baseirrel} and Proposition
\ref{prodpot} imply the following formula:
\begin{equation}\label{GijklPDE}
\begin{split}
\calG(i j \mid k l) =
&\sum_{s} \left( \frac{\partial^3 \calN}{\partial y_i \partial y_j
\partial y_s}
         \frac{\partial^3 \calP}{\partial y_k \partial y_l \partial
y_{s^*}^{*}}
         +\frac{\partial^3 \calN}{\partial y_k \partial y_l \partial
y_s}
         \frac{\partial^3 \calP}{\partial y_i \partial y_j \partial
y_{s^*}^{*}}
         \right) \\
& + 18 \sum_{s,t} \frac{\partial^3 \calN}{\partial y_i \partial y_j
\partial y_s}
                \frac{\partial^2 \calT}{\partial y_{s^*}^{*}\partial
y_{t^*}^{*}}
                \frac{\partial^3 \calN}{\partial y_k \partial y_l
\partial y_t}.
\end{split}
\end{equation}
Thus the equation of Theorem \ref{Gijklprop} is a third-order partial
differential equation.

\section{Determining the characteristic numbers}
\label{recform}

\par
We now use Theorem \ref{Gijklprop} and formula (\ref{GijklPDE}),
applied
with $Y_i = Y_j = h$, $Y_k = Y_l = h^2$ to give a recursive scheme
for determining the Gromov-Witten invariants
$\langle (h^2)^{3d-3} \gamma_1\rangle_d$ and
$\langle (h^2)^{3d-3} \gamma_2 \cdot \gamma_3\rangle_d$ where
$\gamma_1 \in A^3(S_2)$ and $\gamma_2, \gamma_3 \in A^2(S_2)$.  To
shorten expressions, we will employ the following conventions:  If
$\calM$ is a potential, $\calM_i$ will denote the partial derivative
$\partial \calM/ \partial y_i$ where $y_i$ is a coefficient with
respect to the $z$-basis, whereas $\calM_{i^*}$ will denote
$\partial \calM/ \partial y_{i^*}^{*}$ where $y_{i^*}^{*}$ is a
coefficient with respect to the $i$-basis. For example,
(\ref{GijklPDE}) abbreviates to
\begin{equation}\label{abbrev}
\calG(i j \mid k l) =
   \sum_{s} \left( \calN_{ijs} \calP_{kls^*} + \calN_{kls} \calP_{ijs^*}
         \right)
  + 18 \sum_{s,t} \calN_{ijs} \calT_{s^* t^*} \calN_{klt}.
\end{equation}
\par First we note that with $Y_i = Y_j = h$, $Y_k = Y_l = h^2$, the
classical potential $\calP$ contributes just a single term to the
first sum in (\ref{abbrev}), namely a $1$ that arises when $i = j =
100$ and that corresponds to the triple product
$$
h \cdot h \cdot \hd i= Y_{100}\cdot Y_{100} \cdot Y_{011}^{*}
$$
in $A^*(S_2)$.  Thus Theorem \ref{Gijklprop} yields
$$
\calN_{(200)^3} =
   18 \sum_{s,t} \left(
           \calN_{100,200,s}\calT_{s^* t^*}\calN_{100,200,t}
             - \calN_{100,100,s}\calT_{s^* t^*}\calN_{200,200,t}
                 \right).
$$
Expressed in terms of $N^{(d)}$ and $E^{(d)}$, this gives a formula in
each degree:
\begin{equation*}
\begin{split}
N^{(d)}_{(200)^3}E^{(d)} =
  18 \sum_{\substack{d_1+d_2 = d\\ d_1,d_2 > 0\\ s,t}}
     &\left\{
        \left( N^{(d_1)}E^{(d_1)}\right)_{100,200,s}
        \calT_{s^* t^*}
        \left( N^{(d_2)}E^{(d_2)}\right)_{100,200,t} \right. \\
   &\left. -  \left( N^{(d_1)}E^{(d_1)}\right)_{100,100,s}
        \calT_{s^* t^*}
        \left( N^{(d_2)}E^{(d_2)}\right)_{200,200,t}
     \right\}.
\end{split}
\end{equation*}
We simplify to obtain
\begin{equation}\label{recPDE}
\begin{split}
N^{(d)}_{(200)^{3}} =
  \frac{18}{E^{(d)}} \sum_{\substack{d_1+d_2 = d\\ d_1,d_2 > 0\\ s,t}}
     &\left\{
     d_1 d_2 \left( N^{(d_1)}_{200}E^{(d_1)}\right)_s
             \calT_{s^* t^*}
             \left( N^{(d_2)}_{200}E^{(d_2)}\right)_t \right. \\
    &\left. -\: d_1^2 \left( N^{(d_1)}E^{(d_1)}\right)_s
            \calT_{s^* t^*}
            \left( N^{(d_2)}_{200,200}E^{(d_2)}\right)_t
     \right\}.
\end{split}
\end{equation}

\begin{thm}\label{recprop}
Given $N^{(1)}$, equation (\ref{recPDE}) uniquely determines
$N^{(d)}_{(200)^{3d-3}}$ for all $d >0$.
\end{thm}

\begin{proof}
Apply $\partial^{3d-6}/ \partial y_{200}^{3d-6}$ to equation
(\ref{recPDE}). From the first term on the right we obtain terms of
the form
$$
d_1d_2\frac{\partial^{a_1}}{\partial y_{200}^{a_1}}
      \left( N^{(d_1)}_{200}\right) E^{(d_1)} \calT_{s^* t^*}
      \frac{\partial^{a_2}}{\partial y_{200}^{a_2}}
      \left( N^{(d_2)}_{200}\right) E^{(d_2)}
$$
and from the second term we obtain terms of the form
$$
d_1^2\frac{\partial^{b_1}}{\partial y_{200}^{b_1}}
     \left( N^{(d_1)}\right) E^{(d_1)} \calT_{s^* t^*}
     \frac{\partial^{b_2}}{\partial y_{200}^{b_2}}
     \left( N^{(d_2)}_{200,200}\right) E^{(d_2)},
$$
where $a_1 + a_2 = b_1 + b_2 = 3d-6$.  Since
$$
\frac{\partial^{3d}N^{(d)}}{\partial y_{200}^{3d}} = 0,
$$
any nonzero terms must have $a_i \leq 3d_i -2$ for $i = 1,2$.  Thus
$a_1 = 3d-6-a_2 \geq 3d_1 -4$ and $a_2 \geq 3d_2 -4$.  Similarly,
$b_1 \leq 3d_1 - 1$, so $b_2 = 3d-6-b_1 \geq 3d_2 - 5$ and
$b_2 + 2 \leq 3d_2 - 1$, so $b_1 = 3d-6-b_2 \geq 3d_1-3$.  To
summarize, $N^{(d)}_{(200)^{3d-3}}$ is
determined by $N^{(d')}_{(200)^{a}}$
for $d' < d$
and $a \geq 3d' -3$.  The theorem follows by induction.
\end{proof}

\begin{cor}\label{GWcor}
Let $\gamma_1 \in A^3(S_2)$ and $\gamma_2, \gamma_3 \in A^2(S_2)$.
Then equation (\ref{recPDE}) determines the Gromov-Witten invariants
$\langle (h^2)^{3d-3} \gamma_1\rangle_d$ and
$\langle (h^2)^{3d-3} \gamma_2 \cdot \gamma_3\rangle_d$
for all $d >0$, given the invariants in the case $d=1$.
\end{cor}

\begin{proof}
Recall that
$$
N^{(d)} = \sum_{\|\mathbf{a}\| = 3d-1}
          N_d(\mathbf{a}) \frac{\mathbf{y}^{\mathbf{a}}}{\mathbf{a}!}
$$
where $N_d(\mathbf{a}) = \langle \prod_i \gamma_i\rangle_d$  (See the
beginning of {\S}\ref{potpde} for notation.)  Thus
$$
N^{(d)}_{(200)^{3d-3}} = \sum_{\|\mathbf{b}\| = 2}
          N_d(\mathbf{a}) \frac{\mathbf{y}^{\mathbf{b}}}{\mathbf{b}!}
$$
where $\mathbf{b} = \mathbf{a} - (3d-3,0,\ldots,0)$.  Hence
if $a_{200} \geq 3d-3$
we obtain
the Gromov-Witten invariant $N_d(\mathbf{a})$ by computing
$\mathbf{b}!$ times the coefficient of
$\mathbf{y}^{\mathbf{b}}$ in $N^{(d)}_{(200)^{3d-3}}$.
\end{proof}

\par
The Gromov-Witten invariants in degree~$1$ are
\begin{gather*}
  \langle h^2\hd\rangle_1 = 1;             \quad \!
  \langle h^2z\rangle_1 = 3;               \quad \!
  \langle \hd^2z \rangle_1 = -3;        \\
  \langle h^2\cdot h^2\rangle_1 = 1;       \quad \!
  \langle h^2\cdot \hd^2\rangle_1 = 0;     \quad \!
  \langle h^2\cdot hz\rangle_1 = 0;        \quad \!
  \langle h^2\cdot \hd z\rangle_1 = -3;    \quad \!
  \langle \hd^2\cdot \hd^2\rangle_1 = 0; \\
  \langle \hd^2\cdot hz\rangle_1 = 0;      \quad \!
  \langle \hd^2\cdot \hd z\rangle_1 = 0;   \quad \!
  \langle hz\cdot hz\rangle_1 = 0;         \quad \!
  \langle hz\cdot \hd z\rangle_1 = 0;      \quad \!
  \langle \hd z\cdot \hd z\rangle_1 = 9.
\end{gather*}

\noindent
These numbers are easily obtained by replacing
$z$ by $i-3h+3\hd$ and noting that the lift of a line to $S_2$
does not meet the divisor at infinity $I$. Thus
$$
N^{(1)} = y_{210} + 3y_{201} - 3y_{021} + \frac{1}{2}y_{200}^2 -
3y_{200}y_{011} + \frac{9}{2}y_{011}^2.
$$
\par
Using (\ref{recPDE}), we have obtained the corresponding thirteen
Gromov-Witten invariants through degree~$6$; the values are shown in
Table \ref{GWtable}. The characteristic numbers
that do not involve $z$ were
also obtained in \cite{ErnstromKennedy}.
There are also several cases overlapping those of
Caporaso and Harris \cite{CH}
and we have checked that our values agree with theirs.
\par
Note that in five instances the entries of a row of
Table \ref{GWtable} are three times
as large as those in a preceding row; for example
\begin{equation*}
\langle (h^2)^{3d-3} \cdot hz \cdot \hd z \rangle_d
= 3
\langle (h^2)^{3d-3} \cdot \hd^2 \cdot \hd z \rangle_d.
\end{equation*}
This is easily explained by using $hz-3\hd^2=hi$
and then employing an \textit{ad hoc} argument to show that
\begin{equation*}
\langle (h^2)^{3d-3} \cdot hi \cdot \hd z \rangle_d
= 0.
\end{equation*}
Here is a rough outline of the argument.
The class is represented by a cycle
supported on two loci, the first of which is
the locus of curves which are unions of curves
of lower degree, passing through $3d-3$ specified points,
meeting at a point on a specified line, and having a flex
tangent line through a specified point.
The other is the locus of curves having a cusp on a
specified line,
passing through $3d-3$ points, and having a flex
tangent line through a specified point. For general data, these loci are empty.
Similar arguments apply in the other four instances.

\begin{table}
\begin{center}
{\renewcommand{\arraystretch}{1.25}
\begin{tabular}{ll||cccccc}
  & & \multicolumn{6}{c}{$d$} \\
  & & 1 & 2 & 3 & 4 & 5 & 6 \\ \hline \hline
$\gamma_1$
 & $h^2\hd$ & 1 & 1 & 10 & 428 & 51040 & 13300176 \\
 & $h^2z$   & 3 & 3 & 30 & 1284 & 153120 & 39900528 \\
 & $\hd^2z$ & $-3$ & 0 & 21 & 1452 & 216180 & 64150200 \\ \hline
$\gamma_2 \cdot \gamma_3$
 & $h^2\cdot h^2$ & 1 & 1 & 12 & 620 & 87304 & 26312976 \\
 & $h^2\cdot \hd^2$ & 0 & 2 & 36 & 2184 & 335792 & 106976160 \\
 & $h^2\cdot hz$ & 0 & 6 & 108 & 6552 & 1007376 & 320928480 \\
 & $h^2\cdot \hd z$ & $-3$ & 0 & 54 & 4872 & 894528 & 315755712 \\
 & $\hd^2\cdot \hd^2$ & 0 & 4 & 100 & 7200 & 1222192 & 415085088 \\
 & $\hd^2\cdot hz$ & 0 & 12 & 300 & 21600 & 3666576 & 1245255264 \\
 & $\hd^2\cdot \hd z$ & 0 & 0 & 150 & 15912 & 3223944 & 1214002800 \\
 & $hz\cdot hz$ & 0 & 36 & 900 & 64800 & 10999728 & 3735765792 \\
 & $hz\cdot \hd z$ & 0 & 0 & 450 & 47736 & 9671832 & 3642008400 \\
 & $\hd z\cdot \hd z$ & 9 & 0 & 63 & 22860 & 6556140 & 2948122440
 \medskip
\end{tabular}}
\end{center}
\caption{The Gromov-Witten invariants
$\langle (h^2)^{3d-3} \gamma_1\rangle_d$ and
$\langle (h^2)^{3d-3} \gamma_2 \cdot\gamma_3\rangle_d$ where
$\gamma_1 \in A^3(S_2)$ and $\gamma_2, \gamma_3 \in A^2(S_2)$ for $1
\leq d \leq 6$.}
\label{GWtable}
\end{table}

\section{Contact formulas}
\label{contactform}

\par
In this final section, we present two enumerative formulas for
simultaneous triple contacts between fixed plane curves and members
of a family of plane curves.  The first formula is a generalization of
Theorem 4 of \cite{CKsimult} to arbitrarily many fixed curves.
(In that theorem we said we were using the fiber product of
two copies of the second-order stable lift of the family of curves,
but in fact an examination of the proof shows that we were
using the join, i.e., the unique component dominating the fiber
product of two copies of the family of curves.)
\par
Before we state the result, we note that the action of the projective
general linear group $PGL(2)$ on $\pp$ lifts to compatible actions on
$S_1$ and $S_2$.  This is well-known for $S_1$,
the incidence correspondence for $\pp$.  Moreover, the action of
$PGL(2)$ on $S_1$ respects fibers and hence the induced action on
$\mathbf{P}(TS_1)$ must send a focal plane to another focal plane.
Thus $PGL(2)$ acts naturally on $S_2$.

\par
To set notation, let $C_1,\ldots,C_n$ denote fixed curves in $\pp$
and $(C_1)_2,\ldots,(C_n)_2$ their lifts to $S_2$.  Suppose $\calX$
is a family of plane curves over $T$, a complete parameter space of dimension
$2n$. As in
{\S}\ref{GWinv}, let $\jxnt$
denote the join over $T$ of $n$ copies of $\calX_2$.
Let $(S_2)^n$ be the product of $n$ copies of the Semple bundle
variety.
Let $\pi_i\colon  (S_2)^n \to S_2$ be projection onto
the $i$th factor and let $\sigma\colon  (S_2)^n \times T \to (S_2)^n$ also
denote projection.

\begin{thm}\label{triplecontactprop}
Suppose that the reduced plane curves $C_1,\ldots,C_n$
each contain no line and that
the general member of $\calX$ is reduced and  contains no
line.  If
$C_1,\ldots,C_n$ are in general position with respect to
the action of $(PGL(2))^n$ on $(S_2)^n$, then the number of
simultaneous triple contacts between $C_1,\ldots,C_n$, and
members of
$\calX$ is given by
$$
\int_{(S_2)^n} \pi_1^{*}[(C_1)_2] \cup \cdots \cup
\pi_n^{*}[(C_n)_2]
               \cap \sigma_{*}[\jxnt].
$$
Equivalently, by the projection formula, this number is
$$
\int_{(S_2)^n \times T}
\sigma^{*}\pi_1^{*}[(C_1)_2] \cup \cdots \cup
\sigma^{*}\pi_n^{*}[(C_n)_2]
               \cap [\jxnt].
$$
\end{thm}
\begin{proof}
It was shown in \cite{RS2} and \cite[Theorem 1]{CKtrip}
that $S_2$ has three orbits under the action of $PGL(2)$: a dense
orbit $\calO(-)$, represented by the germ at the origin of a
nonsingular curve without a flex; a three-dimensional orbit
$\calO(0) = Z$, represented by the germ at the origin of a line; and
another three-dimensional orbit $\calO(\infty) = I$, represented by
the germ at the origin of a curve with an ordinary cusp.  Thus
$(S_2)^n$ has $3^n$ orbits of the form $\calO = \calO(a_1) \times
\cdots \times \calO(a_n)$ where each $a_i$ can be one of the
symbols $-,0,\infty$.

\par
We will show that for every nondense orbit $\calO$ on $(S_2)^n$
\begin{equation}\label{orbitineq}
\dim \left( (C_1)_2 \times \cdots \times (C_n)_2 \cap \calO\right)
   + \dim \left( \jxnt \cap (\calO \times T)\right)
   < \dim \calO.
\end{equation}
Hence, since $C_1,\ldots,C_n$ are assumed to be in
general position with respect to the action of $(PGL(2))^n$,
the transversality theory of \cite{Kleiman4} shows that
$$
(C_1)_2 \times \cdots \times (C_n)_2
   \cap \sigma(\jxnt) \cap \calO = \emptyset.
$$
Thus $(C_1)_2 \times \cdots \times (C_n)_2$ and
$\sigma(\jxnt)$ intersect transversely and all intersections
must occur in the dense orbit $\calO(-) \times \cdots \times
\calO(-)$.  The product of lifts
$(C_1)_2 \times \cdots \times (C_n)_2$ is the closure of the graph of
function defined on a dense subset of $C_1 \times \cdots \times C_n$;
the join $\jxnt$ is likewise the closure of the graph of a
function defined on a dense subset of $\calX_n(T)$, the fiber product
over $T$ of $n$ copies of $\calX$.  Thus, by general
position, all intersections between
$(C_1)_2 \times \cdots \times (C_n)_2$ and $\sigma(\jxnt)$
are intersections between the graphs.  Therefore every intersection
point is an $n$-tuple $(x_1,\ldots,x_n)$ in which $x_1$
lies over a nonsingular point of $C_1$ and some member $X_t$ of $\calX$,
$x_2$ lies over a nonsingular point of $C_2$ and $X_t$, etc.  Since each
$x_i$ is a point of $S_2$, we see that the second-order data of $C_i$ and
$X_t$ at $x_i$ must be identical for $i = 1,\ldots,n$.  Thus each
intersection point is a simultaneous triple contact.

\par
Establishing (\ref{orbitineq}) is quite easy.  First, if
$\calO = \calO(a_1) \times \cdots \times \calO(a_n)$, then
$$
\dim \calO = 4n-q
$$
where $q$ is the number of $a_i$'s which are not the
symbol ``$-$''.  Next, since none of the curves $C_i$ contains a
line, we have
$$
\dim\left( (C_1)_2 \times \cdots \times (C_1)_2 \cap \calO \right)
    \leq r
$$
where $r$ is the number of $a_i$'s which are ``$-$''.  Finally, since
the general member of $\calX$ is reduced and contains no line,
$$
\dim\left(\jxnt \cap(\calO \times T) \right)
< \dim\left(\jxnt\right) = 3n,
$$
from which (\ref{orbitineq}) follows.
\end{proof}

There is a natural generalization of Theorem~\ref{triplecontactprop}
allowing a mixture of single, double and triple contact conditions. The proof
is by a transversality argument just as above and is therefore omitted.

\begin{thm}\label{mixedcontactprop}
Let $P_1,\dots,P_r$ be points in $\pp$; let $C_1,\dots,C_s$
be any reduced plane curves; let $D_1,\dots,D_t$ be reduced
curves containing no lines. Assume that all these data are in
general position with respect to the action of $(PGL(2))^n$ on
$(\pp)^n$, where $n=r+s+t$. Let $\calX$ be a complete family
of plane curves whose general member is reduced and contains
no lines, over an $m$-dimensional parameter space where $m=r+s+2t$.
Then the number of members of
$\calX$ passing through the points $P_i$, tangent
to the curves $C_i$ and making a triple contact with each of
the curves $D_i$ is equal to
$$
\int_{(S_2)^n \times T}
\prod_{i=1}^r \sigma^* \pi_i^*(f_2^*f_1^*[P_i])
\prod_{i=1}^s \sigma^* \pi_{r+i}^*(f_2^*[(C_i)_1])
\prod_{i=1}^t \sigma^* \pi_{r+s+i}^*[(D_i)_2] \cap [\jxnt].
$$
\end{thm}

\par
We now apply Theorem~\ref{mixedcontactprop} to the
(completed) family ${\overline \calR}(d)$ of all
rational plane curves of degree~$d$. Using
Theorem~\ref{comparison}, the number of members of this
family satisfying the conditions of
Theorem~\ref{mixedcontactprop}, with $m=3d-1$,
is equal to the Gromov-Witten invariant
\begin{equation}\label{mixedGW}
\left\langle\prod_{i=1}^r f_2^*f_1^*[P_i]
\cdot\prod_{i=1}^s f_2^*[(C_i)_1]
\cdot\prod_{i=1}^t [(D_i)_2]\right\rangle_d.
\end{equation}
 From \cite[p.\ 182]{CKsimult} we recall the identities
$$
[(C_i)_1]=d(C_i)\check{h}^2+\check{d}(C_i)h^2
$$
and
$$
[(D_i)_2]=d(D_i)\check{h}^2z+\check{d}(D_i)h^2z+\kappa(D_i)h^2\check{h}
$$
where $d(C)=\int h\cap [C]$ is the degree, $\check{d}(C)=\int \check{h} \cap
[(C)_1]$ is the class, and
$\kappa(C)=\int i\cap [(C)_2]$ is the number of cusps (assuming no
worse singularities than ordinary nodes and cusps).
Plugging these identities into (\ref{mixedGW}), we get an expression
in terms of the degree~$d$ Gromov-Witten invariants of the family, and
the three invariants of each
specified curve.
\par
In particular, we obtain an explicit formula for the number $N_d(C)$
of rational curves of degree $d$ passing though $3d-3$ points and
making a triple contact with one specified curve $C$ (assumed to be reduced
and containing no lines).
Strictly speaking, the formula applies to the completed family
${\overline \calR}(d)$. But, as we now argue, all of the curves
counted by the formula are honest rational curves; there are no
contributions from reducible curves. Indeed, suppose that
$d=d_1+d_2$. If $d_1>1$, then according to Theorem~\ref{mixedcontactprop}
there are only finitely members of ${\overline \calR}(d_1)$ making a triple
contact with $C$ and passing through $3d_1-3$ points, and there are only
finitely many members of ${\overline \calR}(d_2)$ passing through
$3d_2-1$ points. Since $(3d_1-3)+(3d_2-1)<3d-3$, there is no union
of two such curves satisfying all the specified conditions. Similarly,
there are finitely many flex tangent lines to $C$, and finitely many
members of ${\overline \calR}(d-1)$ passing through $3d-4$ points.
Thus there is no union of a line and another curve satisfying all the
specified conditions.
\par
To avoid a clash of notations we use $c$ and $\check{c}$ for
the degree and class of $C$. Then
\begin{equation}\label{final}
N_d(C) =
c\left\langle(h^2)^{3d-3}\cdot\check{h}^2z\right\rangle
+\check{c}\left\langle(h^2)^{3d-3}\cdot h^2z\right\rangle
+\kappa\left\langle(h^2)^{3d-3}\cdot h^2\check{h}\right\rangle.
\end{equation}
Using inputs from Table~\ref{GWtable}, we obtain the explicit formulas
shown in Table~\ref{triplecontacttable}. (The formula for $d=1$, which
can also be derived from the classical Pl{\"u}cker formulas,
counts the number of flexes on $C$.)

\begin{table}
\begin{center}
{\renewcommand{\arraystretch}{1.25}
\begin{tabular}{c|c}
 $d$ & $N_d(C)$ \\ \hline
1 & $-3c+3\check{c}+\kappa$ \\
2 & $3\check{c}+\kappa$ \\
3 & $21c+30\check{c}+10\kappa$ \\
4 & $1452c+1284\check{c}+428\kappa$ \\
5 & $216180c+153120\check{c}+51040\kappa$ \\
6 & $64150200c+39900528\check{c}+13300176\kappa$
\medskip
\end{tabular}}
\end{center}
\caption{The number of rational plane curves of degree~$d$ passing through
$3d-3$ points and making a triple contact with a curve $C$.}
\label{triplecontacttable}
\end{table}

\providecommand{\bysame}{\leavevmode\hbox to3em{\hrulefill}\thinspace}

\end{document}